\documentclass{cmslatex}
\usepackage[paperwidth=7in, paperheight=10in, margin=.875in]{geometry}
 \usepackage[backref,colorlinks,linkcolor=red,anchorcolor=green,citecolor=blue]{hyperref}
\usepackage{amsfonts,amssymb}
\usepackage{amsmath}
\usepackage{graphicx}
\usepackage{cite}
\usepackage{enumerate}
\usepackage{algorithm,algorithmic}
\graphicspath{{./pics/}}
\usepackage{dsfont}
\usepackage{booktabs}
\renewcommand{\theequation}{\arabic{section}.\arabic{equation}}

\newcommand{\R}{\mathbb{R}}

\allowdisplaybreaks
\providecommand{\norm}[1]{\left\lVert#1\right\rVert}
\renewcommand{\d}{{\,\rm d}}
\providecommand{\abs}[1]{\left\lvert#1\right\rvert}

\begin{document}
\title{Recovery of a Distributed Order Fractional Derivative in an Unknown Medium\thanks{Received date, and accepted date (The correct dates will be entered by the editor).}}
\author{Bangti Jin\thanks{Department of Mathematics, The Chinese University of Hong Kong, Shatin, New Territories, Hong Kong, P.R. China (\texttt{bangti.jin@gmail.com, b.jin@cuhk.edu.hk}) The work of B. J. is supported by UK EPSRC grant EP/T000864/1
and EP/V026259/1, and a start-up fund of The Chinese University of Hong Kong.}
          \and Yavar Kian\thanks{Aix Marseille Universit\'{e}, Universit\'{e} de Toulon, CNRS, CPT, Marseille, France (\texttt{yavar.kian@univ-amu.fr}) The work of Yavar Kian was supported by  the French National Research Agency ANR (project MultiOnde) grant ANR-17-CE40-0029.}}

\pagestyle{myheadings} \markboth{WEIGHT RECOVERY}{Bangti Jin and Yavar Kian} \maketitle

\begin{abstract}
In this work, we study an inverse problem of recovering information about the weight
in distributed-order time-fractional diffusion from the observation at one single point
on the domain boundary. In the absence of an explicit knowledge of the medium, we
prove  that the one-point observation can uniquely determine the support bound of the weight.
The proof is based on asymptotics of the data, analytic continuation and
Titchmarch convolution theorem. When the medium is known, we give an
alternative proof of an existing result, i.e., the one-point boundary observation uniquely
determines the weight. Several numerical experiments are also presented to complement the analysis.\\
\end{abstract}


\section{Introduction}
This work is concerned with an inverse problem for distributed-order time-fractional diffusion.
Let $\Omega\subset \mathbb{R}^d$ ($d\ge 2$) be an open bounded and connected subset with a
$C^{2\lceil\frac{d}{4}\rceil+2}$ boundary $\partial \Omega$ (the notation $\lceil r\rceil$ denotes
the smallest integer $\geq r\in\mathbb{R}$). Let $a:=(a_{ij})_{i,j=1}^d \in
C^{1+2\lceil\frac{d}{4}\rceil} (\overline{\Omega};\mathbb{R}^{d\times d})$ be symmetric
and fulfill the standard ellipticity condition, i.e., there exists a $c>0$ such that
\begin{equation}\label{ell}
   \sum_{i,j=1}^d a_{i,j}(x) \xi_i \xi_j \geq c |\xi|^2, \quad
  \forall x \in \overline{\Omega},\ \xi=(\xi_1,\ldots,\xi_d) \in \R^d,
\end{equation}
and $q \in C^{2\lceil\frac{d}{4}\rceil}(\overline{\Omega})$ is strictly positive in
$\overline{\Omega}$. Then we define the operator $\mathcal{A}$ by
\begin{equation*}
  \mathcal{A} u(x) :=-\sum_{i,j=1}^d \partial_{x_i}
   \left( a_{i,j}(x) \partial_{x_j} u(x) \right)+q(x)u(x),\  x\in\Omega.
\end{equation*}
Finally, $\mu\in C[0,1]$ is a non-negative function satisfying supp$(\mu)=[b_1,b_2]$
with $0\leq b_1<b_2<1$ and $b_1$ and $b_2$ being lower and upper bounds of the
support of the weight $\mu$. We define the function
\begin{equation*}
   K_\mu(t):=\int_0^1\frac{t^{-\alpha}}{\Gamma(1-\alpha)}\mu(\alpha)\d\alpha,
\end{equation*}
where $\Gamma(z)=\int_0^\infty s^{z-1}e^{-s}\d s$, for $\Re(z)>0$, denotes Euler's Gamma function.
Then we define an integral operator $I_{K_\mu}$ by
\begin{equation*}
  I_{K_\mu}g(x,t)=\int_0^tK_\mu(t-s)g(x,s)\d s,\quad g\in L^1_{loc}(\R_+;L^2(\Omega)),\ x\in\Omega,\ t\in\R_+:=(0,+\infty).
\end{equation*}
We define Riemann-Liouville and Djrbashian-Caputo fractional derivatives with the kernel $K_\mu$ for $g\in C([0,+\infty);L^2(\Omega))$ by
\begin{equation*}
  D^{[\mu]}_tg(x,t)=\partial_tI_{K_\mu}g(x,t),\quad \partial^{[\mu]}_tg(x,t)=D^{[\mu]}_t[g-g(\cdot,0)](x,t),\quad \ x\in\Omega,\ t\in\R_+.
\end{equation*}
For $g\in W^{1,1}_{loc}([0,+\infty);L^2(\Omega))$, $\partial^{[\mu]}_tg=I_{K_\mu}
\partial_tg$, which coincides with the classical Djrbashian-Caputo fractional derivative
for the kernel $K_\mu$. Now we consider the following initial boundary value problem:
\begin{equation}\label{eq1}
\left\{ \begin{aligned}
(\partial^{[\mu]}_t+\mathcal{A} ) u(x,t) & =  F(x,t), && \mbox{in }\Omega\times\R_+,\\
 \mathcal R u(x,t) & =  0, && \mbox{on } \partial\Omega\times\R_+, \\
 u(x,0) & =  u_0(x), && \mbox{in } \Omega.
\end{aligned}
\right.
\end{equation}
Here $\mathcal R$ denotes either Dirichlet trace $\mathcal Ru=u_{|\partial
\Omega\times\R_+}$ or conormal derivative $\partial_{\nu_a}$ associated with $a$, i.e.,
$\mathcal Ru=\partial_{\nu_a}u_{|\partial\Omega\times\R_+}=\sum_{i,j=1}^da_{ij}\partial_{x_j}u\nu_i|_{\partial\Omega\times\R_+}$,
where $\nu\in\mathbb{R}^d$ denotes  the unit outward normal to
$\partial\Omega$. Throughout, the adjoint trace $\mathcal{R}^*$ is the Dirichlet trace
if $\mathcal R$ is the Neumann one and vice versa.

Physically, problem \eqref{eq1} describes ultra-slow anomalous diffusion processes, in which the mean
square displacement admits a logarithmic growth \cite{MainardiMura:2008,MeerschaertScheffler:2006}.
The model is also called distributed-order fractional diffusion, generalizing
time-fractional diffusion \cite{P,KubicaRyszewskaYamamoto:2020,Jin:2021book}.
The distributed-order fractional derivative with a constant weight was used by Caputo
\cite{Caputo:1969} in 1969 to describe the stress-strain relation of anelastic media. Later it was used to
formally determine eigenfunctions of torsional modes of anelastic or dielectric spherical
shells and infinite plates and to estimate split periods of their free modes \cite{Caputo:1989}.
The related differential equation was first solved by Caputo \cite{Caputo:1995} using Laplace
transform. Since then, the model has received much attention
in applied disciplines, e.g., polymer physics and kinetics of particles moving in quenched random force
fields \cite{ChechkinGorenfloSokolov:2003,Naber:2004,SokolovChechkinKlafter:2004}.
There is also a growing body of literature \cite{Kochubei:2008,Luchko:2009fcaa,Kochubei:2009jp,
LiLuchkoYamamoto:2014fcaa,LiLuchkoYamamoto:2017,KubicaRyszewska:2018,LIKianSoccorsi:2019,
KubicaRyszewska:2019} on the analysis of the direct problem. Kochubei \cite{Kochubei:2008,
Kochubei:2009jp} investigated fundamental solutions to the Cauchy problems
for both ODEs and PDEs with $\mu\in C^1[0,1]$. The long- and short-time asymptotic behaviors were discussed in
\cite{LiLuchkoYamamoto:2014fcaa} using an argument similar to Watson lemma; see also
\cite{KubicaRyszewska:2018} for a more general case using energy estimates.

The weight $\mu$ characterizes ultra-slow diffusion processes. In practice, it is determined
by the inhomogeneity of the medium, but the relevant physical law is still unknown at present, and moreover,
it also cannot be measured directly. This naturally motivates the inverse problem of identifying
the weight $\mu$ from, e.g., the measurement of the solution $u$
at one point. When the medium is known and $F\equiv0$, Li et al
\cite{LiLuchkoYamamoto:2017} proved that the measurement at one interior point can determine
the weight  $\mu$. The proof is based on analytic properties of solutions,
especially time analyticity, and properties of Laplace transform of
solutions of \eqref{eq1}.  In the 1D case, Rundell and Zhang \cite{RundellZhang:2017}
studied a similar inverse problem using the value of the
solution $u$ in the time interval $(0,\infty)$, where the main tool in the analysis is a fractional
version of the classical $\theta$ function. More recently, Li et al \cite{LiFujishiroLi:2020} studied the inverse
problem with nonzero boundary conditions, and generalized the result in \cite{RundellZhang:2017}
using a Harnack type inequality and strong maximum principle.
Also  \cite{LiuSunYamamoto:2021} is the only work on numerical recovery,
based on Tikhonov regularization and conjugate gradient method.
See \cite{LiLiuYamamoto:2019} for an overview of
results on determining the orders / weights in time-fractional models.

In this work, we investigate the inverse problem of recovering the support info
of the weight $\mu$, i.e., $b_1$ and $b_2$, from the observation at one point on the boundary
$\partial\Omega$, where the medium is unknown. While the approach in \cite{LiLuchkoYamamoto:2017}
can be successfully applied to the case of a known medium, it is unclear if similar results
still hold when the medium is unknown. The main goal of
this work is to determine the support bounds $b_1$ and $b_2$ of the weight $\mu$ from the
knowledge $\mathcal R^*u(x_0,t)$ for $t\in(0,T)$, $x_0\in\partial\Omega$, when the medium
is unknown. We answer
the question in the affirmative. See Section \ref{sec:main} for the precise statements of
the main results. Further, we provide numerical experiments to
complement the analysis for both known and unknown media. The numerical recovery
is still under-explored, due to the complexity associated with approximating
$\partial_t^{[\mu]}u$. In the appendix, we describe an
easy-to-implement algorithm to perform numerical simulation.

The rest of the paper is organized as follows. In Section \ref{sec:main}, we present
and discuss the main results. Then in Section \ref{sec:prelim}, we
recall preliminary results. In
Section \ref{sec:asymptotic}, we present several asymptotic results for
$t\to 0^+$ and $t\to +\infty$. Then in Section \ref{sec:proof-main} we give the proofs
of the main results. We provide several numerical
experiments in Section \ref{sec:numer}. Finally
we give concluding remarks in Section \ref{sec:concl}.

\section{Main results and discussions} \label{sec:main}
In this section, we describe the main theoretical results. We
define an operator $A=\mathcal A$ acting on $L^2(\Omega)$ with its domain $D(A)$ given by
$D(A)=\{v\in H^2(\Omega):\ \mathcal A v\in L^2(\Omega),\ \mathcal R v=0\mbox{ on }\partial\Omega\}.$
By \cite[Theorem 2.5.1.1]{Grisvard:1985}, for all $\ell=1,\ldots,\lceil\frac{d}{4}\rceil+1$,
\begin{equation}\label{DA}
   D(A^\ell)=\{v\in H^{2\ell}(\Omega):\ \mathcal Rv=\mathcal R(\mathcal Av)=\ldots=\mathcal R(\mathcal A^{\ell-1} v)=0\}.
\end{equation}
The spectrum of $A$ consists of a non-decreasing sequence of strictly
positive eigenvalues $(\lambda_{n})_{n\geq1}$ (repeated with multiplicity). We
introduce an $L^2(\Omega)$ orthonormal basis of eigenfunctions
$(\varphi_{n})_{n\geq1}$ of $A$ associated with the eigenvalues $(\lambda_{n})_{n\geq1}$.
For all $s\geq 0$, we define the operator $A^s$  by
\begin{equation*}
   A^s h=\sum_{n=1}^{+\infty}\left\langle h,\varphi_n\right\rangle \lambda_{n}^s\varphi_{n},\quad h\in D(A^s)
  = \Big\{h\in L^2(\Omega):\ \sum_{n=1}^{+\infty}\abs{\left\langle h, \varphi_{n}\right\rangle}^2 \lambda_{n}^{2s}<\infty\Big\},
\end{equation*}
where $\left\langle \cdot,\cdot\right\rangle$ denotes the $L^2(\Omega)$
inner product, and the graph norm by
$\|h\|_{D(A^s)}= (\sum_{n=1}^{+\infty}\abs{\langle h,
  \varphi_{n}\rangle}^2 \lambda_{n}^{2s})^{\frac{1}{2}}$,  for $h\in D(A^s)$. It
  is convenient to introduce the concept of an admissible tuple, to
describe regularity conditions on the problem data so that the direct problem has a weak solution
$u$ in the sense of Definition \ref{def:weak} below.
\begin{definition}\label{def:admissible}
A tuple $(\Omega,a,q,f,u_0)$ is said to be admissible if the following conditions are fulfilled.
\begin{itemize}
\item[{\rm(i)}] $\Omega\subset \mathbb{R}^d$ is a $C^{2\lceil\frac{d}{4}\rceil+2}$ bounded open set,
   $a:=(a_{i,j})_{ i,j=1}^d \in C^{1+2\lceil\frac{d}{4}\rceil}(\overline{\Omega}; \mathbb{R}^{
   d\times d})$ satisfies the ellipticity condition \eqref{ell}, {and} $q \in C^{2\lceil\frac{d}{4}\rceil}
   (\overline{\Omega})$ is strictly positive on $\overline{\Omega}$.
\item[{\rm(ii)}] $u_0\in D(A^{r+1})$, with $r>\frac{d+3}{2}$, $F\in L^1(\R_+;L^2(\Omega))$ and there
exists $T\in \R_+$ such that
\begin{equation} \label{F}
   F(x,t):=\sigma(t)f(x),\quad (x,t)\in\Omega\times(0,T),
\end{equation}
with $f\in D(A^r)$ and $\sigma\in L^1(0,T)$ and $\sigma\not\equiv0$.
\end{itemize}
\end{definition}

Following \cite{Kian:2021-wellposed}, we define weak solutions of problem \eqref{eq1} as follows.
\begin{definition}\label{def:weak}
Let the coefficients  in \eqref{eq1} satisfy \eqref{ell}. A function
$u\in L_{loc}^1(\R_+;L^2(\Omega))$ is said to be a weak solution to problem \eqref{eq1} if it satisfies
the following conditions.
\begin{enumerate}
\item[{\rm(i)}] The  identity $D_t^\mu[u-u_0](x,t)+\mathcal A u(x,t)=F(x,t),$  $t\in\R_+$,  $x\in\Omega$
holds in the sense of distribution in $\Omega\times\R_+$.
\item[{\rm(ii)}] There hold $I_{K_\mu} [u-u_0]\in W_{loc}^{1,1}(\R_+;H^{-2}(\Omega))$ and $I_{K_\mu}
[u-u_0](x,0)=0$,  $x\in\Omega$.
\item[{\rm(iii)}] There holds $\tau_0=\inf\{\tau>0:\  e^{-\tau t}u\in L^1(\R_+;L^2(\Omega))\}<\infty$ and
there exists $\tau_1\ge \tau_0$ such that for all $p\in\mathbb C$ with $\Re( p)>\tau_1$,
$\widehat u(\cdot,p):=\int_0^\infty e^{-p t}u(\cdot,t)\d t\in H^2(\Omega)$, and $\mathcal R \widehat u(\cdot,p)=0$ on $\partial\Omega$.
\end{enumerate}
\end{definition}

The  conditions in Definition \ref{def:weak} describe different aspects of  \eqref{eq1}:
(i) is associated with the governing equation, (ii) describes
the link with the initial condition and (iii) gives the boundary condition.

The main results can be stated as follows. In formula \eqref{t1ab}, $\delta_{ij}$ denotes the Kronecker
delta symbol.
\begin{assumption}\label{ass:R}
$\mathcal R$ is the Dirichlet trace or $\mathcal R$ is the Neumann trace
and the following condition holds
\begin{equation}\label{t1ab}
  a_{ij}(x)=\delta_{ij},\quad x\in\partial\Omega.
\end{equation}
\end{assumption}
\begin{theorem}\label{thm:main1}
Let $(\Omega_k,a_k,q_k,f_k,0)$, $k=1,2$, be two admissible tuples with
$f_k\not\equiv0$, $\mu_k\in C[0,1]$ non-negative such that
supp$(\mu_k)=[b_1^k,b_2^k]$ with $0\leq b_1^k<b_2^k<1$, and $u^k$ the
corresponding weak solution of problem \eqref{eq1} with $\mu=\mu_k$. If
there exist $x_k\in\partial\Omega_k$, $k=1,2$, such that the condition
\begin{equation}\label{t1a}
   \min\left(|\mathcal R_1^*f_1(x_1)|,|\mathcal R_2^*f_2(x_2)|\right)>0,
\end{equation}
holds. Then,  the condition
\begin{equation}\label{t1b}
  \mathcal R_1^*u^1(x_1,t)=\mathcal R_2^*u^2(x_2,t),\quad t\in(0,T)
\end{equation}
implies $b_2^1=b_2^2$. In addition, if Assumption \ref{ass:R} holds for $\mathcal{R}_k$ and
the following sign condition
\begin{equation}\label{t1c}
   f_1\geq 0 \textrm{ and }f_2\geq 0\quad \textrm{or} \quad f_1\leq 0 \textrm{ and }f_2\leq 0
\end{equation}
holds, then condition \eqref{t1b} implies also $b_1^1=b_1^2$.
\end{theorem}

\begin{theorem}\label{thm:main2}
Let $(\Omega_k,a_k,q_k,0, u_0^k)$, $k=1,2$, be two admissible tuples with
$u_0^k\not\equiv0$, $\mu_k\in C[0,1]$ non-negative such that supp$(\mu_k)=
[b_1^k,b_2^k]$ with $0\leq b_1^k<b_2^k<1$, and $u^k$ the corresponding weak solution of
problem \eqref{eq1} with $\mu=
\mu_k$. If there exist  $x_k\in\partial\Omega_k$, $k=1,2$, such that the condition
\begin{equation}\label{t2a}
\min\left(|\mathcal R_1^*\mathcal A_1u^1_0(x_1)|,|\mathcal R_2^*\mathcal A_2u^2_0(x_2)|\right)>0
\end{equation}
holds, then for any $T_1\in(0,T)$ and $T_2\in (T_1,T)$, with $T$ given by \eqref{F}, the condition
\begin{equation}\label{t2b}
\mathcal R_1^*u^1(x_1,t)=\mathcal R_2^*u^2(x_2,t),\quad t\in(T_1,T_2)
\end{equation}
implies $b_2^1=b_2^2$. In addition, if Assumption \ref{ass:R} holds for $\mathcal{R}_k$, and the following sign condition
\begin{equation}\label{t2c}
   u_0^1\geq 0 \textrm{ and }u_0^2\geq 0\quad \textrm{or} \quad u_0^1\leq 0 \textrm{ and }u_0^2\leq 0
\end{equation}
holds, then condition \eqref{t2b} implies also $b_1^1=b_1^2$.
\end{theorem}

\begin{remark}
The proof in Section \ref{sec:proof-main} indicates that the regularity condition $\mu\in C[0,1]$ can be
relaxed to $\mu$ being piecewise smooth. Notably Theorems \ref{thm:main1} and
\ref{thm:main2} do not require a full knowledge of the problem data $(\Omega,a,q,f,u_0)$.
Uniqueness in the cases of single-term and multi-term models was established in
\cite{JinKian:2021siam} and \cite{JinKian:2021pa}, respectively; see also \cite{Yamamoto:2021} for
a uniqueness result from interior data.
\end{remark}

In the spirit of \cite[Theorem 2.2]{LiLuchkoYamamoto:2017}, we can strengthen the results
when the medium is known.
Then we can show the unique recovery of the weight $\mu\in L^1(0,1)$.
\begin{theorem}\label{thm:known}
Let $(\Omega,a,q,f,u_0)$ be an admissible tuple, with $T$ given by
condition \eqref{F}, $\mu\in L^1(0,1)$ a non-negative
function satisfying
\begin{equation} \label{mumu}
\exists \alpha_0 \in(0,1),\ \exists \varepsilon \in (0,\alpha_0),\ \forall \alpha \in (\alpha_0-\varepsilon,\alpha_0),\quad \mu(\alpha) \ge \frac{\mu(\alpha_0)}{2}>0.
\end{equation}
with $\mu=\mu_k$, and $u^k$ the corresponding weak solution of problem \eqref{eq1}.
Let Assumption \ref{ass:R} hold for $\mathcal{R}$. Fix any point
$x_0\in\partial\Omega$. If $u_0\equiv0$ and $f\not\equiv0$ is of constant sign,
then the condition
\begin{equation}\label{t3a}
\mathcal R^*u^1(x_0,t)=\mathcal R^*u^2(x_0,t),\quad t\in(0,T)
\end{equation}
implies $\mu_1=\mu_2$. Moreover, if $f\equiv0$ and $u_0\not\equiv0$ is of constant sign, for any
$T_1,T_2\in[0,T]$ with $T_1<T_2$, the condition
\begin{equation}\label{t3aa}
\mathcal R^*u^1(x_0,t)=\mathcal R^*u^2(x_0,t),\quad t\in(T_1,T_2)
\end{equation}
implies $\mu_1=\mu_2$.
\end{theorem}

\begin{remark}
Note that in all the results, the data actually does not access the source $F$
for $\Omega\times(T,+\infty)$. Thus, the statement of the inverse problem is on the set
$\Omega\times(0,T)$ instead of $\Omega\times\R_+$.
\end{remark}

\section{Preliminaries}\label{sec:prelim}

In this section we give preliminary results of problem \eqref{eq1}, e.g., existence
of a weak solution and time analyticity. These results will play a crucial role in the uniqueness
proof in Section \ref{sec:proof-main}. Since $\mu\in C[0,1]$ is non-negative and supp$(\mu)=[b_1,b_2]$, for
all $\alpha_0\in[b_1,b_2)$ and all $\varepsilon\in(0,b_2-\alpha_0]$, there exist $\delta_0
\in(0,\varepsilon)$ and $\delta_1\in (0,\delta_0)$ such that
\begin{equation}\label{mu}
   \mu(\alpha_1)>\frac{\mu(\alpha_0+\delta_0)}{2}>0,\quad \forall \alpha_1\in [\alpha_0+\delta_1,\alpha_0+\delta_0].
\end{equation}
Throughout, we fix $\theta\in(\frac{\pi}{2},\min(\frac{\pi}{2b_2},\pi))$, $\delta\in\R_+$ and the
contour $\gamma(\delta,\theta)\subset \mathbb C$ defined by
$\gamma(\delta,\theta):=\{\delta e^{{\rm i}\beta}:\ \beta\in[-\theta,\theta]\} \cup \{r e^{\pm {\rm i}\theta}:\  r\geq\delta\}$,
oriented in the counterclockwise direction. Let $\theta_1\in(0,\frac{2\theta-\pi}{8})$. For any $p\in\mathbb C\setminus(-\infty,0]$ we set $p^\alpha:=e^{\alpha \log(p)}$ with $\log$ being complex logarithm defined on  $\mathbb C\setminus(-\infty,0]$. For all $z\in
D_{\theta_1}:=\{re^{{\rm i}\beta}:\ r>0,\ \beta\in (-\theta_1,\theta_1)\}$, we can define two solution
operators $S_1(z), S_2(z)\in\mathcal B(L^2(\Omega))$ by
\begin{align}
    S_1(z)v&=\frac{1}{2{\rm i}\pi}\int_{\gamma(\delta,\theta_1)} e^{z p}\left(A+\mathcal V(p)\right)^{-1}p^{-1}\mathcal V(p)v \d p,\label{S1}\\
    S_2(z)v&=\frac{1}{2{\rm i}\pi}\int_{\gamma(\delta,\theta_1)}e^{z p}\left(A+\mathcal V(p)\right)^{-1}v \d p,\label{S}
\end{align}
with the weight function $\mathcal{V}(p)$ given by
$\mathcal V(p)=\int_0^1p^\alpha\mu(\alpha)\d\alpha$, for $p\in\mathbb C\setminus(-\infty,0]$.
The operators $S_1(z)$ and $S_2(z)$ correspond to the initial data and right
hand side, respectively.

Under condition \eqref{mu}, the argument of  \cite[Theorem 1.4]{LIKianSoccorsi:2019} and \cite[Lemma 3.1]{Kian:2021-wellposed} gives the following result.
\begin{lemma}\label{l1}
For all $\tau\in[0,1)$, the map $z\mapsto S_j(z)$, $j=1,2$, is holomorphic in the set $D_{\theta_1}$
as a map taking values in $\mathcal B(L^2(\Omega); D(A^\tau))$ and there exist $C>0$ and $\delta'\in(0,b_2
-b_1)$ such that
\begin{align*}
\norm{S_1(z)}_{B(L^2(\Omega); D(A^\tau))}&\leq C\max(|z|^{(b_1+\delta')(1-\tau)-1},|z|^{b_2(1-\tau)-1}),\quad z\in D_{\theta_1},\\
\norm{S_2(z)}_{B(L^2(\Omega); D(A^\tau))}&\leq C\max(|z|^{(b_1+\delta')(1-\tau)-b_2},|z|^{b_2(1-\tau)-b_1}),\quad z\in D_{\theta_1}.
\end{align*}
\end{lemma}

Lemma \ref{l1} and \cite[Theorem 1.4]{Kian:2021-wellposed} imply that for $F\in L^1
(\mathbb{R}_+;L^2(\Omega))$ and $u_0\in L^2(\Omega)$,  problem \eqref{eq1} admits a unique weak solution
$u\in L^1_{loc}(\R_+;D(A^\tau))$, $\tau\in [0,1)$, given by
\begin{equation}\label{sol}
u(\cdot,t)=S_1(t)u_0+\int_0^t S_2(t-s)F(\cdot,s) \d s,\quad t\in\R_+.
\end{equation}
Similar to \cite{LIKianSoccorsi:2019}, for $z\in D_{\theta_1}$ and $j=1,2$, we have
\begin{equation*}
   S_j(z)h=\sum_{n=1}^{+\infty}S_{j,n}(z)\left\langle h,\varphi_n\right\rangle\varphi_n,
\end{equation*}
with
$S_{1,n}(z)=\frac{1}{2{\rm i}\pi}\int_{\gamma(\delta,\theta)} e^{z p}(\lambda_n+\mathcal V(p))^{-1}
    p^{-1} \mathcal V(p) \d p$
and $S_{2,n}(z)=\frac{1}{2{\rm i}\pi}\int_{\gamma(\delta,\theta)}e^{z p}(\lambda_n+\mathcal V(p))^{-1} \d p.$
The argument for Lemma \ref{l1} implies also that for all $n\in\mathbb N$
and $j=1,2$, $S_{j,n}$ is holomorphic on $D_{\theta_1}$. Moreover,
 we have  for all $\tau\in[0,1]$ \cite[Lemmas 2.1 and 2.2]{LIKianSoccorsi:2019}
\begin{align}
\abs{S_{1,n}(z)}&\leq C\lambda_n^{-\tau}\max(|z|^{(b_1+\delta')(1-\tau)-1},|z|^{b_2(1-\tau)-1}),\quad z\in D_{\theta_1},\ n\in\mathbb N,\label{l1d}\\
\abs{S_{2,n}(z)}&\leq C\lambda_n^{-\tau}\max(|z|^{(b_1+\delta')(1-\tau)-b_2},|z|^{b_2(1-\tau)-b_1}),\quad z\in D_{\theta_1},\ n\in\mathbb N,\label{l1b}
\end{align}
with the constant $C>0$ depending only on $\theta$, $\mathcal A$, $\mu$  and $\Omega$.
Then we can prove the following representation of the measurement $\mathcal{R}^*u(x,t)$.

\begin{lemma}\label{l2}
Let $f,u_0\in D(A^{\frac{d}{4}})$, $\sigma\in L^1(0,T)$ and $F\in L^1(\R_+;L^2(\Omega))$ satisfy condition
\eqref{F}. Then the map{\color{blue}s} $t\mapsto S_1(t)u_0$ and $t\mapsto S_2(t)f$  are analytic with respect to $t\in\R_+$
as a function taking values in $C^1(\overline{\Omega})$. Moreover, problem \eqref{eq1} admits a unique weak
solution $u\in L^1_{loc}(\R_+; L^2(\Omega))$ such that the restriction of $u$ to $\Omega\times (0,T)$ belongs to  $L^1(0,T;C^1(\overline{\Omega}))$ and
\begin{equation}\label{l2a}
\mathcal R^* u(x,t)=\mathcal R^* [S_1(t)u_0](x)+ \int_0^t\sigma(s) \mathcal R^* [S_2(t-s)f](x)\d s,\quad t\in(0,T),\ x\in\partial\Omega.
\end{equation}
\end{lemma}
\begin{proof}
We only prove \eqref{l2a} for $u_0\equiv0$. By the identity
\eqref{DA} and interpolation, the space $D(A^{\frac{d}{4}+\frac{3}{4}})$ embeds continuously
into $H^{\frac{d}{2}+\frac{3}{2}}(\Omega)$ and Sobolev embedding theorem implies that $D(A^{\frac{d}{4}
+\frac{3}{4}})$ embeds continuously into $ C^1(\overline{\Omega})$. In addition, by applying \eqref{l1b},
for all  $z\in D_{\theta_1}$ and all $m_1,m_2\in\mathbb N$, $m_1<m_2$, we have
\begin{align}\label{l2b}
\norm{\sum_{n=m_1}^{m_2}S_{2,n}(z)\left\langle f,\varphi_n\right\rangle\varphi_n}_{ C^1(\overline{\Omega})}&\leq C\norm{\sum_{n=m_1}^{m_2}S_{2,n}(z)\left\langle f,\varphi_n\right\rangle\varphi_n}_{D(A^{\frac{d}{4}+\frac{3}{4}})}\\
 &\leq C\max(|z|^{\frac{b_1+\delta'}{4}-1},|z|^{\frac{b_2}{4}-1})\left(\sum_{n=m_1}^{m_2}\lambda_n^{\frac{d}{2}}|\left\langle f,\varphi_n\right\rangle|^2\right)^{\frac{1}{2}},\nonumber
\end{align}
with $C>0$ independent of $z$, $m_1$ and $m_2$. This and the condition $f\in D(A^{\frac{d}{4}})$
yield that the sequence
$\sum_{n=1}^{N}S_{2,n}(z)\left\langle f,\varphi_n\right\rangle\varphi_n$, $N\in\mathbb{N}$,
converges uniformly with respect to $z$ on any compact subset of $D_{\theta_1}$ to $S_2(z)f$ as a function taking
values in $C^1(\overline{\Omega})$. Thus, the map $D_{\theta_1}\ni z\mapsto S_2(z)f$ is holomorphic
as a function taking values in $C^1(\overline{\Omega})$. This shows the first assertion.
Next, in view of \eqref{F} and \eqref{sol}, we have
\begin{equation*}
   u(\cdot,t)=\int_0^t S_2(t-s)F(\cdot,s) \d s=\int_0^tS_2(t-s)\sigma(s)f\d s=\int_0^t\sigma(s)S_2(t-s)f\d s,\quad t\in(0,T).
\end{equation*}
Thus, the estimate \eqref{l2b} implies
\begin{equation*}
   \norm{u(\cdot,t)}_{D(A^{\frac{d}{4}+\frac{3}{4}})}\leq C\norm{f}_{D(A^{\frac{d}{4}})}\left(\max(t^{\frac{b_1+\delta}{4}-1},t^{\frac{b_2}{4}-1})
   \mathds{1}_{(0,T)}\right)*(|\sigma|\mathds{1}_{(0,T)})(t),\quad t\in(0,T),
\end{equation*}
where $\mathds{1}_{(0,T)}$ denotes the characteristic function of $(0,T)$ and $*$ denotes
the convolution product. Hence, Young's inequality gives $u\in L^1(0,T;
D(A^{\frac{d}{4}+\frac{3}{4}}))\subset L^1(0,T; C^1(\overline{\Omega}))$, showing
the second assertion. In the same way, by applying  \eqref{l2b},
we deduce the representation \eqref{l2a}.
\end{proof}

Below we also use the following result from \cite[Lemma 2.3]{JinKian:2021pa}.
\begin{lemma}\label{l3}
The following estimate and identity hold
\begin{align}
    \sum_{n=1}^{+\infty}\lambda_n|\left\langle v,\varphi_n\right\rangle|\norm{\mathcal R^*\varphi_n}_{L^\infty(\partial\Omega)}\leq C\norm{v}_{D(A^{\frac{d}{2}+\frac{3}{2}})},\quad v\in D(A^{\frac{d}{2}+\frac{3}{2}}), \label{l3a}\\
    \mathcal R^*v(x)=\sum_{n=1}^{+\infty} \left\langle v,\varphi_n\right\rangle \mathcal R^*\varphi_n(x),\quad  x\in\partial\Omega,\quad v\in D(A^{\frac{d}{2}+\frac{1}{2}}).\label{l3b}
\end{align}
\end{lemma}

\section{Asymptotic properties of solutions as $t\to0$ and $t\to+\infty$}\label{sec:asymptotic}
In this section, we fix an admissible tuple $(\Omega,a,q,h,u_0)$ with  $\mu\in C[0,1]$ a non-negative function
satisfying supp$(\mu)=[b_1,b_2]$, $0\leq b_1<b_2<1$. Let $w
(\cdot,t)=S_2(t)h\in L^1_{loc}(\R_+;C^1(\overline{\Omega}))\cap C(\R_+;C^1(
\overline{\Omega}))$, cf. Lemma \ref{l2}. Below
we fix $x_0\in\partial\Omega$ and study the asymptotics of $\mathcal R^*w(x_0,t)$
as $t\to0^+$ and $t\to+\infty$. The next result gives the asymptotic behavior of
$\mathcal R^*w(x_0,t)$ as $t\to0^+$.
\begin{lemma}\label{l5}
Let $Q(t,p)$ and $P(t)$ be defined, for $t\in\R_+$, $p\in\mathbb C\setminus(-\infty,0]$, {respectively} by
\begin{equation}\label{l5a}
   Q(t,p):=\int_{0}^1t^{-\alpha}p^\alpha \mu(\alpha)\d\alpha\quad\mbox{and}\quad P(t)=\int_{b_1}^{b_2}t^{-\alpha}\mu(\alpha)\d\alpha.
\end{equation}
Define also
\begin{equation} \label{l5h}
   \mathfrak{J}(t):=\frac{1}{2{\rm i}\pi}\int_{\gamma(1,\theta)}e^{p}Q(t,p)^{-1}\d p.
\end{equation}
Then for $h\in D(A^{\frac{d}{2}+\frac{3}{2}})$, the following asymptotic expansion holds
\begin{equation} \label{l5b}
   \mathcal R^*w(x_0,t)=t^{-1}\mathfrak{J}(t)\mathcal R^*h(x_0)+\underset{t\to0^+}{O}\left(t^{-1}P(t)^{-2}\right).
\end{equation}
\end{lemma}
\begin{proof}
For any $h\in D(A^{\frac{d}{2}+\frac{1}{2}})$, we have $\left(A+\mathcal V(p)\right)^{-1}h
\in D(A^{\frac{d}{2}+\frac{1}{2}})$ and Lemma \ref{l3} implies
\begin{equation*}
   \mathcal R^*\left(A+\mathcal V(p)\right)^{-1}h(x_0)=\sum_{n=1}^{+\infty} \frac{\left\langle h,\varphi_n\right\rangle}{\lambda_n+\mathcal V(p)}\mathcal R^*\varphi_n(x_0).
\end{equation*}
Then \cite[Lemmas 2.1 and 2.2]{LIKianSoccorsi:2019} and Lebesgue dominated
convergence theorem \cite[Theorem 1.19, p. 28]{Evans:2015} imply
\begin{equation*}
     \mathcal R^*w(x_0,t)=\mathcal R^*S_2(t)h(x_0)= \frac{1}{2{\rm i}\pi}\int_{\gamma(\delta,\theta)}e^{t p} \left(\sum_{n=1}^{+\infty} \frac{\left\langle h,\varphi_n\right\rangle}{\lambda_n+\mathcal V(p)}\mathcal R^*\varphi_n(x_0)\right)\d p,
\end{equation*}
where $\delta>0$ is arbitrary. By fixing $\delta=t^{-1}$ and changing variables, we deduce
\begin{align*}
    \mathcal R^*w(x_0,t)&= \frac{1}{2{\rm i}\pi}\int_{\gamma(t^{-1},\theta)}e^{t p} \left(\sum_{n=1}^{+\infty} \frac{\left\langle h,\varphi_n\right\rangle}{\lambda_n+\mathcal V(p)}\mathcal R^*\varphi_n(x_0)\right)\d p\\
     \ &= \frac{1}{2{\rm i}\pi}\int_{\gamma(1,\theta)}e^{p} \left(\sum_{n=1}^{+\infty}\frac{t^{-1}\left\langle h,\varphi_n\right\rangle\mathcal R^*\varphi_n(x_0)}{\lambda_n+Q(t,p)}\right)\d p.
\end{align*}
Next we fix $n\in\mathbb N$, $t\in (0,1)$, $p\in\gamma(1,\theta)$.
Direct computation gives
\begin{equation*}
\begin{aligned}
&\frac{t^{-1}\left\langle h,\varphi_n\right\rangle\mathcal R^*\varphi_n(x_0)}{\lambda_n+Q(t,p)}=\frac{t^{-1}Q(t,p)^{-1}\left\langle h,\varphi_n\right\rangle\mathcal R^*\varphi_n(x_0)}{1+\lambda_nQ(t,p)^{-1}}\\
=&t^{-1}Q(t,p)^{-1}\left\langle h,\varphi_n\right\rangle\mathcal R^*\varphi_n(x_0)-\int_0^1\frac{t^{-1}Q(t,p)^{-2}\lambda_n\left\langle h,\varphi_n\right\rangle\mathcal R^*\varphi_n(x_0)}{(1+s\lambda_nQ(t,p)^{-1})^2}\d s.
\end{aligned}
\end{equation*}
Now summing this expression over $n\in\mathbb N$ and applying Lemma \ref{l3} lead to
\begin{align*}
   &\sum_{n=1}^{+\infty}\frac{t^{-1}\left\langle h,\varphi_n\right\rangle\mathcal R^*\varphi_n(x_0)}{\lambda_n+Q(t,p)}\\
   =&t^{-1}Q(t,p)^{-1}\sum_{n=1}^{+\infty}\left\langle h,\varphi_n\right\rangle\mathcal R^*\varphi_n(x_0)-\sum_{n=1}^{+\infty}\int_0^1\frac{t^{-1}Q(t,p)^{-2}\lambda_n\left\langle h,\varphi_n\right\rangle\mathcal R^*\varphi_n(x_0)}{(1+s\lambda_nQ(t,p)^{-1})^2}\d s\\
   =&t^{-1}Q(t,p)^{-1}\mathcal R^*h(x_0)-\sum_{n=1}^{+\infty}\int_0^1\frac{t^{-1}Q(t,p)^{-2}\lambda_n\left\langle h,\varphi_n\right\rangle\mathcal R^*\varphi_n(x_0)}{(1+s\lambda_nQ(t,p)^{-1})^2}\d s.
\end{align*}
Since $\theta\in(\frac{\pi}{2},\min(\frac{\pi}{2b_2},\pi))$ and supp$(\mu)\subset[b_1,b_2]$,
for all $s\in[0,1]$, $r\in[1,+\infty)$ and all $\beta\in(-\theta,\theta)$, we have
\begin{align*}
  |1+s\lambda_nQ(t,re^{{\rm i}\beta})^{-1}|&=|Q(t,re^{{\rm i}\beta})|^{-1}|Q(t,re^{{\rm i}\beta})+s\lambda_n|\\
  &\geq|Q(t,re^{{\rm i}\beta})|^{-1}|s\lambda_n+\Re(Q(t,re^{{\rm i}\beta}))|\\
  &\geq|Q(t,re^{{\rm i}\beta})|^{-1}\left(s\lambda_n+\int_0^1t^{-\alpha}r^\alpha\cos(\alpha\beta)\mu(\alpha)\d\alpha\right)\\
  &\geq |Q(t,re^{{\rm i}\beta})|^{-1}\cos(b_2\theta) \int_{b_1}^{b_2}t^{-\alpha}r^\alpha\mu(\alpha)\d\alpha\\
  &\geq \cos(b_2\theta)|Q(t,re^{{\rm i}\beta})|^{-1}|Q(t,re^{{\rm i}\beta})|=\cos(b_2\theta).
\end{align*}
Therefore, applying again Lemma \ref{l3} and the above estimate lead to
\begin{align*}
   &\abs{\sum_{n=1}^{+\infty}\frac{t^{-1}\left\langle h,\varphi_n\right\rangle\mathcal R^*\varphi_n(x_0)}{\lambda_n+Q(t,p)}-t^{-1}Q(t,p)^{-1}\mathcal R^*h(x_0)}\\
   &=\abs{\sum_{n=1}^{+\infty}\int_0^1\frac{t^{-1}Q(t,p)^{-2}\lambda_n\left\langle h,\varphi_n\right\rangle\mathcal R^*\varphi_n(x_0)}{(1+s\lambda_nQ(t,p)^{-1})^2}\d s}\\
   &\leq t^{-1}\cos(b_2\theta)^{-2}|Q(t,p)|^{-2}\sum_{n=1}^{+\infty} \lambda_n\abs{\left\langle h,\varphi_n\right\rangle\mathcal R^*\varphi_n(x_0)}\\
   &\leq Ct^{-1}|Q(t,p)|^{-2}\norm{h}_{D(A^{\frac{d+3}{2}})},
\end{align*}
with the constant $C>0$ depending only on $A$, $\theta$ and $\mu$. This and
the definition of $\gamma(\delta,\theta)$ imply
\begin{align*}
      &\abs{\mathcal R^*w(x_0,t)- \frac{\mathcal R^*h(x_0)}{2{\rm i}\pi}\int_{\gamma(1,\theta)}e^{p}t^{-1}Q(t,p)^{-1}\d p}\\
      &\leq Ct^{-1}\norm{h}_{D(A^{\frac{d+3}{2}})} \left(\int_{-\theta}^{\theta}e^{\cos(\beta)}|Q(t,e^{{\rm i}\beta})|^{-2}\d\beta+\int_1^{+\infty} e^{r\cos(\theta)}|Q(t,re^{\pm {\rm i}\theta})|^{-2}\d r\right)\\
      &:= Ct^{-1}\norm{h}_{D(A^{\frac{d+3}{2}})}({\rm I}_1(t)+{\rm I}_2(t)).
\end{align*}
Next we bound the two terms ${\rm I}_1(t)$ and ${\rm I}_2(t)$. Since $\theta\in(\frac{\pi}{2},\min(
\frac{\pi}{2b_2},\pi))$ and supp$(\mu)\subset[b_1,b_2]$, we deduce that, for all $r\in[1,+\infty)$ and all $\beta\in(-\theta,\theta)$,
\begin{align} 
      |Q(t,re^{{\rm i}\beta})| &\geq \big|\Re (Q(t,re^{{\rm i}\beta}))\big|=\int_{b_1}^{b_2}t^{-\alpha}r^{\alpha}\cos(\alpha\beta)\mu(\alpha)\d\alpha\nonumber\\
      &\geq \cos(b_2\theta)r^{b_2}\int_{b_1}^{b_2}t^{-\alpha}\mu(\alpha)\d\alpha.\label{l5e}
\end{align}
Then it follows that
${\rm I}_1(t)\leq C(\int_{b_1}^{b_2}t^{-\alpha}\mu(\alpha)\d\alpha)^{-2}$ and ${\rm I}_2(t)\leq C(\int_{b_1}^{b_2}t^{-\alpha}\mu(\alpha)\d\alpha)^{-2}$.
Hence
\begin{equation} \label{l5g}
   \begin{aligned}
      &\abs{\mathcal R^*w(x_0,t)- \frac{\mathcal R^*h(x_0)}{2{\rm i}\pi}\int_{\gamma(1,\theta)}e^{p}t^{-1}Q(t,p)^{-1}\d p}\\
      \leq& Ct^{-1}\norm{h}_{D(A^{\frac{d+3}{2}})}\left(\int_{b_1}^{b_2}t^{-\alpha}\mu(\alpha)\d\alpha\right)^{-2}.
   \end{aligned}
\end{equation}
This proves the estimate \eqref{l5b}, and completes the proof of the lemma.
\end{proof}

Next we consider the asymptotic behavior of $\mathcal R^*w(x_0,t)$ as $t\to+\infty$.
\begin{lemma}\label{l6}
Let $Q(t,p)$ and $P(t)$ be defined by \eqref{l5a}, and let
\begin{equation} \label{l6a}
 \mathfrak{K}(t)=\frac{1}{2{\rm i}\pi}\int_{\gamma(1,\theta)}e^{p}Q(t,p)\d p.
\end{equation}
Then for $h\in D(A^{\frac{d}{2}+\frac{1}{2}})$, the following asymptotic expansion holds
\begin{equation} \label{l6b}
   \mathcal R^*w(x_0,t)=-t^{-1}\mathcal R^*A^{-2}h(x_0)\mathfrak{K}(t)+\underset{t\to+\infty}{O}\left(t^{-1}P(t)^2\right).
\end{equation}
\end{lemma}
\begin{proof}
Fix $n\in\mathbb N$, $t\in (1,+\infty)$, $p\in\gamma(1,\theta)$. Using Taylor formula, we find
\begin{align}
    \frac{t^{-1}\left\langle h,\varphi_n\right\rangle\mathcal R^*\varphi_n(x_0)}{\lambda_n+Q(t,p)}
       &=t^{-1}\lambda_n^{-1}\left\langle h,\varphi_n\right\rangle\mathcal R^*\varphi_n(x_0)-t^{-1}\lambda_n^{-2}Q(t,p)\left\langle h,\varphi_n\right\rangle\mathcal R^*\varphi_n(x_0)\nonumber\\
       &\quad +2\int_0^1\frac{t^{-1}Q(t,p)^2\left\langle h,\varphi_n\right\rangle\mathcal R^*\varphi_n(x_0)(1-s)}{(\lambda_n+sQ(t,p))^3}\d s.\label{l5i}
\end{align}
Moreover, for all $s\in[0,1]$, we have
\begin{equation*}
   |\lambda_n+sQ(t,re^{{\rm i}\beta})|\geq \lambda_n+s\int_0^1\cos(\alpha\beta)r^\alpha t^{-\alpha}\mu(\alpha)\d\alpha,\quad r\in[1,+\infty),\ \beta\in(-\theta,\theta),
\end{equation*}
and since supp$(\mu)\subset [b_1,b_2]$ and $\theta<\frac{\pi}{2b_2}$, we obtain
\begin{align*}
       & |\lambda_n+sQ(t,re^{i\beta})|\geq \lambda_n+s\int_{b_1}^{b_2}\cos(\alpha\beta)r^\alpha t^{-\alpha}\mu(\alpha)\d\alpha\\
      \geq &\lambda_n+s\cos(b_2\theta)\int_{b_1}^{b_2}r^\alpha t^{-\alpha}\mu(\alpha)\d\alpha
      \geq \lambda_n,\quad r\in[1,+\infty),\ \beta\in(-\theta,\theta).
\end{align*}
Now summing \eqref{l5i} over $n\in\mathbb N$ and applying Lemma \ref{l3} yield
\begin{align*}
     &\quad\abs{\sum_{n=1}^{+\infty}\frac{t^{-1}\left\langle h,\varphi_n\right\rangle\mathcal R^*\varphi_n(x_0)}{\lambda_n+Q(t,p)}-t^{-1}\mathcal R^*A^{-1}h(x_0)+Q(t,p)t^{-1}\mathcal R^*A^{-2}h(x_0)}\\
    &\leq \sum_{n=1}^{+\infty}\abs{\frac{t^{-1}\left\langle h,\varphi_n\right\rangle\mathcal R^*\varphi_n(x_0)}{\lambda_n+Q(t,p)}-\frac{\left\langle h,\varphi_n\right\rangle}{t\lambda_n}\mathcal R^*\varphi_n(x_0)+\frac{Q(t,p)}{t\lambda_n^2}\left\langle h,\varphi_n\right\rangle\mathcal R^*\varphi_n(x_0)}\\
    &\leq Ct^{-1}|Q(t,p)|^2\norm{h}_{D(A^{\frac{d+1}{2}})}\leq C|p|^{2b_2}t^{-1}\left(\int_{b_1}^{b_2}t^{-\alpha}\mu(\alpha)\d\alpha\right)^2\norm{h}_{D(A^{\frac{d+1}{2}})}.
\end{align*}
Integrating the expression in $p\in\gamma(1,\theta)$ gives
\begin{align*}
&\quad \abs{\mathcal R^*w(x_0,t)- \frac{\mathcal R^*A^{-1}h(x_0)}{t}\frac{1}{2{\rm i}\pi}\int_{\gamma(1,\theta)}e^{p}\d p+\frac{\mathcal R^*A^{-2}h(x_0)}{t}\frac{1}{2{\rm i}\pi}\int_{\gamma(1,\theta)}e^{p}Q(t,p)\d p}\\
&\leq Ct^{-1}\norm{h}_{D(A^{\frac{d+1}{2}})}\left(\int_{-\theta}^{\theta}e^{\cos(\beta)}|Q(t,e^{{\rm i}\beta})|^{2}\d\beta+\int_1^{+\infty} e^{r\cos(\theta)}|Q(t,re^{\pm {\rm i}\theta})|^{2}\d r\right)\\
&\leq Ct^{-1}\norm{h}_{D(A^{\frac{d+1}{2}})}\left(\int_{b_1}^{b_2} t^{-\alpha}\mu(\alpha)\d\alpha\right)^2\left(1+\int_1^{+\infty} e^{r\cos(\theta)}r^{2b_2}\d r\right)\\
 &\leq Ct^{-1}\norm{h}_{D(A^{\frac{d+1}{2}})}\left(\int_{b_1}^{b_2} t^{-\alpha}\mu(\alpha)\d\alpha\right)^2.
\end{align*}
Meanwhile, the Cauchy formula gives $\frac{1}{2{\rm i}\pi}\int_{\gamma(1,\theta)}e^{p}\d p=0$,
and upon fixing $\mathfrak{K}(t)$ as \eqref{l6a}, we get
$$\abs{\mathcal R^*w(x_0,t)+t^{-1}\mathcal R^*A^{-2}h(x_0)\mathfrak{K}(t)}\leq Ct^{-1}\norm{h}_{D(A^{\frac{d+1}{2}})}\left(\int_{b_1}^{b_2} t^{-\alpha}\mu(\alpha)\d\alpha\right)^2,$$
and thus we obtain the estimate \eqref{l6b}.
\end{proof}

We will also need the following two intermediate results.
\begin{lemma}\label{l7}
Let $\mathfrak{J}(t)$ and $\mathfrak{K}(t)$ be given by \eqref{l5h} and \eqref{l6a}, respectively.
Then, there exists a constant $C>0$ depending only on $\mu$ and $\theta$ such that
\begin{align}
\label{l7a} C^{-1}P(t)^{-1}\leq\abs{\mathfrak{J}(t)}\leq CP(t)^{-1},\quad t\in(0,1),\\
\label{l7b} C^{-1}P(t)\leq\abs{\mathfrak{K}(t)}\leq CP(t),\quad t\in(1,+\infty).
\end{align}
\end{lemma}
\begin{proof}
It follows from the proof of Lemmas \ref{l5} and \ref{l6} that the following estimate holds
$\abs{{\mathfrak{J}}(t)}\leq CP(t)^{-1}$, for $t\in(0,1)$. Moreover,  we have
\begin{align*}
\abs{\mathfrak{K}(t)}&\leq C\left(\int_{-\theta}^{\theta}e^{\cos(\beta)}|Q(t,e^{{\rm i}\beta})|\d\beta+\int_1^{+\infty} e^{r\cos(\theta)}|Q(t,re^{\pm {\rm i}\theta})|\d r\right)\\
&\leq C\left(\int_{b_1}^{b_2} t^{-\alpha}\mu(\alpha)\d\alpha\right)\left(1+\int_1^{+\infty} e^{r\cos(\theta)}r^{b_2}\d r\right)\leq C\int_{b_1}^{b_2} t^{-\alpha}\mu(\alpha)\d\alpha,
\end{align*}
which shows the following estimate
\begin{equation*}
   \abs{\mathfrak{K}(t)}\leq C P(t), \quad t\in (1,+\infty). 
\end{equation*}
So we only need to prove the following two estimates
\begin{align*}
C^{-1}P(t)^{-1}\leq\abs{\mathfrak{J}(t)},\quad t\in(0,1)\quad \mbox{and}\quad
C^{-1}P(t)\leq\abs{\mathfrak{K}(t)},\quad t\in(1,+\infty).
\end{align*}
Note that
\begin{equation*} 
|Q(t,p)|\leq |p|^{b_2}\int_{b_1}^{b_2}t^{-\alpha}\mu(\alpha)\d\alpha,\quad p\in\gamma(1,\theta),
\end{equation*}
and, by repeating the arguments for Lemmas \ref{l5} and \ref{l6}, we can
prove the first estimate. In the same way, we obtain the second estimate from \eqref{l5e}.
\end{proof}

\begin{remark}
Note that the argument in this section allows deriving asymptotics of the solutions to
the distributed order model. Indeed, applying Lemmas \ref{l6} and \ref{l7}, we have
$$\abs{\mathcal R^*w(x_0,t)}\leq Ct^{-1}P(t)\norm{h}_{D(A^r)},\quad t>1,$$
with $C>0$ independent of $t$ and $h$.
In addition, we find
$$P(t)\leq \norm{\mu}_{L^\infty(0,1)}\int_{b_1}^{b_2}t^{-\alpha}\d\alpha\leq \norm{\mu}_{L^\infty(0,1)}\left(\frac{t^{-b_1}}{\ln(t)}-\frac{t^{-b_2}}{\ln(t)}\right),\quad t>1$$
and consequently
$$\abs{\mathcal R^*w(x_0,t)}\leq C\frac{t^{-1-b_1}}{\ln(t)}\norm{h}_{D(A^r)},\quad t>1.$$
Similarly, by fixing $\epsilon\in(0,(b_2-b_1)/2)$, we obtain
$$P(t)\geq \left(\inf_{[b_1+\epsilon,b_2-\epsilon]}\mu\right)\int_{b_1+\epsilon}^{b_2-\epsilon}t^{-\alpha}\d\alpha= \left(\inf_{[b_1+\epsilon,b_2-\epsilon]}\mu\right)\left(\frac{t^{-\epsilon-b_1}}{\ln(t)}-\frac{t^{-b_2+\epsilon}}{\ln(t)}\right),\quad t\in(0,1)$$
and Lemmas \ref{l5} and \ref{l7} imply
$$\abs{\mathcal R^*w(x_0,t)}\leq Ct^{-1}P(t)^{-1}\norm{h}_{D(A^r)}\leq Ct^{b_2-1+\epsilon}|\ln(t)|\norm{h}_{D(A^r)},\quad t\in(0,1).$$
Moreover, assuming that $b_2=1$ and $\mu(1)>0$, we obtain
$$\abs{\mathcal R^*w(x_0,t)}\leq  C|\ln(t)|\norm{h}_{D(A^r)},\quad t\in(0,1).$$
These estimates can be compared with that in \cite{LiLuchkoYamamoto:2014fcaa} for the
solution of problem \eqref{eq1} with $F\equiv0$.
\end{remark}

The next result gives the asymptotics of the integral $\int_{b_1}^{b_2}t^{-\alpha}\mu(\alpha) \d\alpha$
as $t\to0^+$ and $t\to+\infty$, which will play a role in determining the support of the weight $\mu$.
\begin{lemma}\label{l8}
The following identities hold
\begin{align}
\label{l8a} \lim_{t\to+\infty} t^b\int_{b_1}^{b_2}t^{-\alpha}\mu(\alpha) \d\alpha=0,\quad b\in(-\infty,b_1),\\
\label{l8b} \lim_{t\to+\infty} t^b\int_{b_1}^{b_2}t^{-\alpha}\mu(\alpha) \d\alpha=+\infty,\quad b\in(b_1,+\infty),\\
\label{l8c} \lim_{t\to0^+} \frac{t^{-b}}{\int_{b_1}^{b_2}t^{-\alpha}\mu(\alpha) \d\alpha}=0,\quad b\in(-\infty,b_2),\\
\label{l8d} \lim_{t\to0^+} \frac{t^{-b}}{\int_{b_1}^{b_2}t^{-\alpha}\mu(\alpha) \d\alpha}=+\infty,\quad b\in(b_2,+\infty).
\end{align}
\end{lemma}
\begin{proof}
To show the identity \eqref{l8a}, we fix $b\in (-\infty,b_1)$. Note that
\begin{equation} \label{l8e}
\int_{b_1}^{b_2}t^{-\alpha}\mu(\alpha) \d\alpha\leq t^{-b_1}\int_{b_1}^{b_2}\mu(\alpha) \d\alpha,\quad t\in(1,+\infty).
\end{equation}
Moreover, since $\mu\in C[0,1]$ is non-negative with supp$(\mu)=[b_1,b_2]$,
$\mu$ cannot be uniformly vanishing and $\int_{b_1}^{b_2}\mu(\alpha) \d\alpha>0.$
Combining this with \eqref{l8e} gives \eqref{l8a}. To show \eqref{l8b},
we fix $b\in (b_1,+\infty)$. By \eqref{mu} with $\alpha_0=b_1$ and $\varepsilon=b-b_1$,
there exist $\delta_0\in(0,\varepsilon)$ and $\delta_1\in (0,\delta_0)$ such that \eqref{mu}
is fulfilled. Then, it follows
\begin{equation*}
\int_{b_1}^{b_2}t^{-\alpha}\mu(\alpha) \d\alpha\geq \int_{b_1+\delta_1}^{b_1+\delta_0}t^{-\alpha}\mu(\alpha) \d\alpha\geq (\delta_0-\delta_1)\frac{\mu(b_1+\delta_0)}{2} t^{-b_1-\delta_1},\quad t\in(1,+\infty).
\end{equation*}
Since $(\delta_0-\delta_1)\frac{\mu(b+\delta_0)}{2}>0$ and $b_1+\delta_1<b_1+
\varepsilon=b$, we deduce easily \eqref{l8b} from this last estimate. For \eqref{l8c}, we fix
$b\in (b_1,b_2)$ and using \eqref{mu} with  $\alpha_0=b$ and $\varepsilon=b_2-b$, we obtain
\begin{equation*}
  \int_{b_1}^{b_2}t^{-\alpha}\mu(\alpha) \d\alpha\geq \int_{b+\delta_1}^{b+\delta_0}t^{-\alpha}\mu(\alpha) \d\alpha\geq (\delta_0-\delta_1)\frac{\mu(b+\delta_0)}{2} t^{-b-\delta_0},\quad t\in(0,1),
\end{equation*}
with $\delta_0\in(0,b_2-b)$ and $\delta_1\in(0,\delta_0)$. This and the inequality
$(\delta_0-\delta_1)\frac{\mu(b+\delta_0)}{2}>0$ and the condition $\delta_0>0$ imply \eqref{l8c}.
Similarly, upon noting $\int_{b_1}^{b_2}t^{-\alpha}\mu(\alpha) \d\alpha\leq t^{-b_2}
\int_{b_1}^{b_2}\mu(\alpha) \d\alpha $ for $t\in(0,1)$, we deduce \eqref{l8d}.
\end{proof}

\section{Proof of Theorems \ref{thm:main1}, \ref{thm:main2} and \ref{thm:known}}\label{sec:proof-main}

To prove Theorems \ref{thm:main1} and \ref{thm:main2}, we need an intermediate result.
\begin{lemma}\label{l9} Let $(\Omega,a,q,f,u_0)$ be an admissible tuple and let the condition
\begin{equation}\label{t3ab}
a_{ij}(x)=\delta_{ij},\quad x\in\partial\Omega,
\end{equation}
 be fulfilled. If $f\not\equiv0$ and $u_0\not\equiv0$ are of constant sign and
 $\mathcal R$ is the Neumann trace. Then, for any $x_0\in\partial\Omega$ we have
 $\mathcal R^*A^{-2}f(x_0)\neq0$ and $\mathcal R^*A^{-1}u_0(x_0)\neq0$.
\end{lemma}
\begin{proof}
We only show the result for  $A^{-2}f$ with $f\leq 0$ since that for $A^{-1}u_0$ follows
 similarly. Let $v=A^{-2}f$. Assume the contradictory, i.e., there exists $x'\in\partial\Omega$ such that
 $v(x')=\mathcal R^*v(x')=0$. Since $f\in D(A^r)$, $v\in D(A^{r+2})$ and,
Sobolev embedding theorem and elliptic regularity theory
\cite[Theorem 2.5.1.1]{Grisvard:1985} imply $v\in C^1(\overline{\Omega})$. The
maximum principle \cite[Theorem 8.1]{GilbargTrudinger:2001} implies
$A^{-1}f\leq0$ and since $f\not\equiv0$, we have $A^{-1}f\not\equiv0$. If $v$ is
constant, then $qv=\mathcal A v=A^{-1}f$ and since $A^{-1}f\not\equiv0$, $A^{-1}f\leq0$
and $q$ is strictly positive, the constant $v$ is negative,
contradicting $v(x')=0$. Hence $v$ is not a constant function. Now
we consider the following two cases.

\noindent\textbf{Case 1.} $\max_{x\in\partial\Omega}v(x)=0$. Since $v$ is
not constant, the strong maximum principle \cite[Theorem 3.5]{GilbargTrudinger:2001} implies
that for all $x\in\Omega$, $v(x)<0=v(x')$ and Hopf's lemma
\cite[Lemma 3.4]{GilbargTrudinger:2001} yields $\partial_\nu v(x')>0$. By
\eqref{t3ab} and the condition $v\in D(A)$, this contradicts the condition $\partial_\nu
v=\mathcal R v=0$ on $\partial\Omega$.

\noindent \textbf{Case 2.} $\max_{x\in\partial\Omega}v(x)\neq0$. Since $v(x')=0$,
we have $\max_{x\in\partial\Omega}v(x)>0$. Since $v\in C(\partial\Omega)$, we
can find $x_2\in\partial\Omega$ such that $v(x_2)=\max_{x\in\partial\Omega}v(x)>0$ and, since
$v$ is not constant, the strong maximum principle implies that for all $x\in\Omega$ we have
$v(x)<v(x_2)$. Then, Hopf's lemma yields $\partial_\nu v(x_2)>0$.
This also contradicts the condition $\partial_\nu v=\mathcal R v=0$ on $\partial\Omega$.

\noindent In both cases we have a contradiction and it follows that for any $x_0\in\partial\Omega$,
 we have $\mathcal R^*A^{-2}f(x_0)\neq0$.
\end{proof}

\subsection{Proof of Theorem \ref{thm:main1}}

We divide the proof into two steps.

\textbf{Step 1.} We prove  that \eqref{t1b} implies $b_2^1=b_2^2$.
By Lemma \ref{l2}, we have
\begin{equation*}
\mathcal R^*_j u_j(x_j,t)=\int_0^t\sigma(s) \mathcal R^*_j [S_2^j(t-s)f_j](x_j)\d s,\quad t\in(0,T),\ j=1,2,
\end{equation*}
with $S_2^j$ corresponding to \eqref{S} with $\mu=\mu_j$ and with $A=A_j$ (acting in
$L^2(\Omega_j)$) and $\mathcal R=\mathcal R_j$. Let $v_j(t)=
\mathcal R^*_j [S_2^j(t)f_j](x_j)$ as a function in $L^1(0,T)$, cf. Lemma \ref{l2}.
Therefore, condition \eqref{t1b} implies
\begin{equation*}
\int_0^t\sigma(s)[v_1(t-s)-v_2(t-s)]\d s=0,\quad t\in(0,T).
\end{equation*}
By Titchmarsh convolution theorem \cite[Theorem VII]{Titchmarsh:1926}, there exist
$t_1,t_2\in[0,T]$ such that $t_1+t_2\geq T$, $\sigma_{|(0,t_1)}\equiv0$ and $(v_1-v_2)_{|(0,t_2)}
\equiv0$. Since $\sigma_{|(0,T)}\not\equiv0$, we have $t_1<T$. This proves
 $t_2=t_1+t_2-t_1\geq T-t_1>0$ and
$\mathcal R^*_1 [S_2^1(t)f_1](x_1)=\mathcal R^*_2 [S_2^2(t)f_2](x_2)$ for $t\in(0,t_2)$.
The analyticity of the map $\R_+\ni t\mapsto \mathcal R^*_j [S_2^j(t)f_j](x_j)$, cf.
Lemma \ref{l2}, implies
$\mathcal R^*_1 [S_2^1(t)f_1](x_1)=\mathcal R^*_2 [S_2^2(t)f_2](x_2)$ for $t\in\R_+$.
Let $w_j(\cdot,t)=S_2^j(t)f_j$. Then
\begin{equation}\label{t1e}
  \mathcal R_1^*w_1(x_1,t)=\mathcal R_2^*w_2(x_2,t),\quad t\in\R_+.
\end{equation}
Now we deduce from Lemma \ref{l5} that
\begin{equation*}
  \mathcal R_j^*w_j(x_j,t)=t^{-1}\mathfrak{J}_j(t)\mathcal R_j^*f_j(x_j)+\underset{t\to0^+}{O}\left(t^{-1}P_j(t)^{-2}\right),
\end{equation*}
with
$Q_j(t,p):=\int_{0}^1t^{-\alpha}p^\alpha \mu_j(\alpha)\d\alpha$, $P_j(t)
   =\int_{b_1^j}^{b_2^j}t^{-\alpha}\mu_j(\alpha)\d\alpha$
and $\mathfrak{J}_j(t):=\frac{1}{2{\rm i}\pi}\int_{\gamma(1,\theta)}e^{p}Q_j(t,p)^{-1}\d p$.
Therefore, \eqref{t1e} implies
\begin{equation}\label{t1f}
t^{-1}\mathfrak{J}_1(t)\mathcal R_1^*f_1(x_1)+\underset{t\to0^+}{O}\left(t^{-1}P_1(t)^{-2}\right)=t^{-1}\mathfrak{J}_2(t)\mathcal R^*_2f_2(x_2)+\underset{t\to0^+}{O}\left(t^{-1}P_2(t)^{-2}\right).
\end{equation}
Combining this with \eqref{t1a}, \eqref{l7a} and \eqref{l8c}-\eqref{l8d} yields $b_2^1=b_2^2$.

\textbf{Step 2.}
Now, we fix $b_2=b_2^1=b_2^2$, and  assume that  \eqref{t1c} holds. Then we prove
that condition \eqref{t1b} implies $b_1^1=b_1^2$. The estimate \eqref{t1e} and Lemma \ref{l6} imply
\begin{align} \label{t1g}
   &-t^{-1}\mathcal R_1^*A_1^{-2}f_1(x_1)\mathfrak{K}_1(t)+\underset{t\to+\infty}{O}\left(t^{-1}P_1(t)^2\right)\\
   =&-t^{-1}\mathcal R_1^*A_2^{-2}f_2(x_2)\mathfrak{K}_2(t)+\underset{t\to+\infty}{O}\left(t^{-1}P_2(t)^2\right),\nonumber
\end{align}
with $\mathfrak{K}_j(t)=\frac{1}{2{\rm i}\pi}\int_{\gamma(1,\theta)}e^{p}Q_j(t,p)\d p.$
By conditions \eqref{l7a} and \eqref{l8a}-\eqref{l8b}, this implies $b_1^1=b_1^2$ if
\begin{equation} \label{t1h}
  \mathcal R_j^*A_j^{-2}f_j(x_j)\neq0,\quad j=1,2.
\end{equation}
So it remains to prove \eqref{t1h}. If $\mathcal R_j$ is the Neumann trace and condition
\eqref{t1ab} holds,  Lemma \ref{l9} implies \eqref{t1h}. Similarly, by
\cite[Lemma 3.1]{JinKian:2021siam}, if  $\mathcal R_j$ is the Dirichlet trace,
condition \eqref{t1c} implies
\begin{equation} \label{t1i}
   \partial_{\nu^j}A_j^{-2}f_j(x_j)\neq0,\quad j=1,2,
\end{equation}
with $\nu^j$ the unit outward normal to $\partial\Omega_j$. Meanwhile, since $\mathcal R_j$ is
the Dirichlet trace, letting the matrix $B_j=(a_{k\ell}^j(x_j))_{ k,
\ell=1}^d$ and since $A_j^{-2}f_j=0$ on $\partial\Omega_j$, we get
\begin{align*}
  &\mathcal R_j^*A_j^{-2}f_j(x_j)=\partial_{\nu_{a_j}} A_j^{-2}f_j(x_j)=[\nabla A_j^{-2}f_j(x_j)]\cdot [B_j\nu^j(x_j)]\\
=&[(\partial_{\nu^j} A_j^{-2}f_j(x_j))\nu^j(x_j)]\cdot [B_j\nu^j(x_j)] =(\partial_{\nu^j} A_j^{-2}f_j(x_j))\nu^j(x_j)^TB_j\nu^j(x_j).
\end{align*}
Moreover, it follows from the ellipticity condition \eqref{ell} that
\begin{equation*}
\nu^j(x_j)^TB_j\nu^j(x_j)= \sum_{k,\ell=1}^d a_{k,\ell}^j(x_j) \nu^j_k(x_j) \nu^j_\ell(x_j) \geq c |\nu^j(x_j)|^2=c>0.
\end{equation*}
This and \eqref{t1i} imply \eqref{t1h}. Thus, in both cases,
\eqref{t1h} holds and we deduce $b_1^1=b_1^2$.

\subsection{Proof of Theorem \ref{thm:main2}}
Note that condition \eqref{F} implies
$u^j(t,\cdot)=S_1^j(t)u_0^j$ for $t\in(0,T)$
with $S_1^j$ corresponding to \eqref{S1} with $\mu=\mu_j$ and with $A=A_j$ (acting in
$L^2(\Omega_j)$) and with $\mathcal R=\mathcal R_j$.
Since $u_0^j\in D(A_j)$, the argument of \cite[pp. 21-22]{Kian:2021-wellposed} implies $u^j\in W^{1,1}(0,T;L^2(\Omega_j))$ with
$\partial_tu^j(t,\cdot)=-S_2^j(t) A_ju_0^j$.
Since $A_ju_0^j\in D(A_j^r)$, $j=1,2$, repeating the argument of Lemma \ref{l2} yields
$u^j\in W^{1,1}(0,T;D(A_j^r))$. Then Lemma \ref{l3} implies
\begin{equation*}
   \partial_t\mathcal R_j^*u^j(t,x_j)=\mathcal R_j^*\partial_tu^j(t,x_j)=-\mathcal R_j^*S_2^j(t) A_ju_0^j(x_j),\quad t\in(0,T).
\end{equation*}
Therefore, differentiating condition \eqref{t1b} in $t$ gives
\begin{equation*}
  \mathcal R_1^*S_2^1(t) A_1u_0^1(x_1)=\mathcal R_2^*S_2^2(t) A_2u_0^2(x_2),\quad t\in (T_1,T_2).
\end{equation*}
This and the time analyticity of the map $t\mapsto \mathcal R_j^*S_2^j(t) A_j
u_0^j(x_j)$, cf. Lemma \ref{l2}, imply \eqref{t1e} with $w_j(\cdot,t)=S_2^j(t)
A_ju_0^j$. Then, repeating the argument for Step 1 of
the proof of Theorem \ref{thm:main1} with $f_j=A_ju_0^j$ yields $b_2^1=b_2^2$. In
addition, if condition \eqref{t2c} holds and $\mathcal R_k$ is
the Dirichlet trace, \cite[Lemma 3.1]{JinKian:2021siam} implies
that for $f_j=A_ju_0^j$, condition \eqref{t1i} holds. Then, similar to
Step 2 of the proof of Theorem \ref{thm:main1}, we deduce \eqref{t1h}. Similarly, if
$\mathcal R_j$ is the Neumann trace and condition \eqref{t1ab} holds,
Lemma \ref{l9} implies  \eqref{t1h}. Hence, in both cases  we deduce from \eqref{t1g} that
$b_1^1=b_1^2$.

\subsection{Proof of Theorem \ref{thm:known}}
For $F\in L^1(\R_+;L^2(\Omega))$, problem \eqref{eq1} admits a unique weak solution $u\in L^1_{loc}(\R_+;L^2(\Omega))$
given by \eqref{sol} \cite[Theorem 1.4]{Kian:2021-wellposed}, and the first statement of Lemma
\ref{l1} and the results of Lemma  \ref{l2} are still valid (\cite[Theorem 1.4]{LIKianSoccorsi:2019} and
\cite[Lemma 3.1]{Kian:2021-wellposed}). Now we can prove Theorem \ref{thm:known}.
\begin{proof}
We prove the result for $u_0\equiv0$ and $f\not\equiv0$ of constant sign. The result for
$f\equiv0$ and $u_0\not\equiv0$ of constant sign follows similarly as  Theorem \ref{thm:main2}.
The argument of Theorem \ref{thm:main1} yields \eqref{t1e} with $x_1=x_2=x_0$ and
$w_j(\cdot,t)=S_2^j(t)f$, $j=1,2$, with $S_2^j$ defined in \eqref{S} with $\mu=\mu_j$.
For all $p\in\mathbb C_+=\{z\in\mathbb C:\ \Re(z)>0\}$, the Laplace
transform $\widehat w_j(\cdot,p)=\int_0^{+\infty} e^{-p t}w_j(\cdot,t)\d t$ of $w_j$
is well defined as an element of $H^2(\Omega)$ (\cite[Proposition 3.2]{Kian:2021-wellposed} and \cite{LIKianSoccorsi:2019}), and we have
\begin{equation} \label{t3b}
   \widehat w_j(\,\cdot\,,p)=\left(A+\mathcal V_j(p)\right)^{-1}f=\sum_{n=1}^{+\infty} \frac{\left\langle f,\varphi_n\right\rangle}{\lambda_n+\mathcal V_j(p)}\varphi_n,\quad\mbox{with }\mathcal V_j(p)=\int_0^1p^\alpha\mu_j(\alpha)\d\alpha.
\end{equation}
Moreover, the argument for Lemma \ref{l2} and \cite[Lemma 3.1]{Kian:2021-wellposed} give
\begin{equation*}
   \norm{w_j(\cdot,t)}_{C^1(\overline{\Omega})}\leq C\norm{w_j(\cdot,t)}_{D(A^r)}\leq
   C\max\left(t^{\alpha_0-\varepsilon-1},t^{\alpha_0-1}\right)\norm{f}_{D(A^r)},\quad t>0,
\end{equation*}
and hence $\widehat{\mathcal R^*w_j(x_0,\cdot)}(p)$ is well defined for $p\in\mathbb C_+$,
and fixing $v_j(t)=\mathcal R^*w_j(x_0,t)$ leads to $\widehat{v_j}(p)=\mathcal R^*\widehat{w_j}(x_0,p)$ for
$p\in\mathbb C_+$. This, \eqref{t3b} and Lemma \ref{l3} imply
\begin{equation*}
   \widehat{v_j}(p)=\sum_{n=1}^{+\infty} \frac{\left\langle f,\varphi_n\right\rangle}{\lambda_n+\mathcal V_j(p)}\mathcal R^*\varphi_n(x_0),\quad p\in\mathbb C_+,\ j=1,2.
\end{equation*}
Therefore, applying the Laplace transform in $t\in\R_+$ to the identity \eqref{t1e} gives
\begin{equation*}
\sum_{n=1}^{+\infty}\frac{\left\langle f,\varphi_n\right\rangle}{\lambda_n+\mathcal V_1(p)}\mathcal R^*\varphi_n(x_0)=\sum_{n=1}^{+\infty} \frac{\left\langle f,\varphi_n\right\rangle}{\lambda_n+\mathcal V_2(p)}\mathcal R^*\varphi_n(x_0),\quad p\in\mathbb C_+,
\end{equation*}
which implies
\begin{equation} \label{t3c}
   (\mathcal V_1(p)-\mathcal V_2(p))\sum_{n=1}^{+\infty} \frac{\left\langle f,\varphi_n\right\rangle}{(\lambda_n+\mathcal V_1(p))(\lambda_n+\mathcal V_2(p))}\mathcal R^*\varphi_n(x_0)=0,\quad p\in\R_+.
\end{equation}
Let
\begin{equation*}
     \mathcal G(p):=\sum_{n=1}^{+\infty}\frac{\left\langle f,\varphi_n\right\rangle}{(\lambda_n+\mathcal V_1(p))(\lambda_n+\mathcal V_2(p))}\mathcal R^*\varphi_n(x_0),\quad p\in\mathbb{R}_+.
\end{equation*}
By Lemma \ref{l3}, $\mathcal G\in C([0+\infty))$ with
\begin{equation*}
   \mathcal G(0)=\sum_{n=1}^{+\infty} \frac{\left\langle f,\varphi_n\right\rangle}{(\lambda_n+\mathcal V_1(0))(\lambda_n+\mathcal V_2(0))}\mathcal R^*\varphi_n(x_0)=\sum_{n=1}^{+\infty} \frac{\left\langle f,\varphi_n\right\rangle}{(\lambda_n)^2}\mathcal R^*\varphi_n(x_0)=\mathcal R^* A^{-2}f(x_0).
\end{equation*}
Then, Lemma \ref{l9}, \cite[Lemma 3.1]{JinKian:2021siam} and the argument at
Step 2 of Theorem \ref{thm:main1} imply $|\mathcal G(0)|>0$ for $\mathcal R$ being
either Dirichlet or  Neumann trace. Hence, there
exists $\varepsilon\in\R_+$ such that $|\mathcal G(p)|>0$, $p\in(0,\varepsilon)$.
This, \eqref{t3c} and fixing $\mu=\mu_1-\mu_2$ yield
\begin{equation} \label{t3d}
  \int_0^1p^\alpha\mu(\alpha)\d\alpha=0,\quad p\in(0,\varepsilon).
\end{equation}
Let $\mu_*$ be the extension of $\mu$ by zero to $\R_+$. Then it follows from \eqref{t3d} that
\begin{equation*}
  \widehat{\mu_*}(-\ln(\tau))=\int_0^{+\infty}e^{\alpha\ln(\tau)}\mu_*(\alpha)\d\alpha=\int_0^1\tau^\alpha
  \mu(\alpha)\d\alpha=0,\quad \tau\in(0,\varepsilon),
\end{equation*}
which implies
\begin{equation} \label{t3e}
\widehat{\mu_*}(p)=0,\quad p\in(-\ln(\varepsilon),+\infty).
\end{equation}
Since $\mu_*\in L^1(\R_+)$ is compactly supported, $\mathbb C\ni p\mapsto \widehat{\mu_*}
(p)$ is holomorphic and \eqref{t3e} implies $\widehat{\mu_*}(p)=0$ for all $p\in\mathbb C$. Then, the
injectivity of the Laplace transform implies $\mu_*\equiv0$, i.e., $\mu_1=\mu_2$.
\end{proof}

\section{Numerical results and discussions}\label{sec:numer}
In this section, we numerically illustrate the analysis,
including both support recovery in Theorems \ref{thm:main1}--\ref{thm:main2} and
weight recovery in Theorem \ref{thm:known}. The exact data $g^\dag=\mathcal{R}^*u(x_0,t)$,
$x_0\in\partial\Omega$, is generated by a fully discrete scheme described in
the appendix. The noisy data $g^\delta$ is generated by $g^\delta(t) =
g^\dag(t) + \epsilon\|g^\dag\|_{L^\infty(0,T)}\xi(t)$, where $\xi(t)$ follows the
standard Gaussian distribution, and $\epsilon\geq0$ denotes the noise level.

\subsection{Bound recovery}
First we illustrate Theorems \ref{thm:main1}
and \ref{thm:main2} with the following example, where
$\mathds{1}_S$ denotes the characteristic function of the set $S$.
\begin{example}\label{exam:bound}
In this example, we take the domain $\Omega=(0,1)$, with
\begin{itemize}
  \item[{\rm(i)}] $u_0(x)=x(1-x)e^x$, $f\equiv0$, and $a(x)=1 + x^2$.
  \item[{\rm(ii)}] $u_0\equiv0$, $f(x,t)= x(1-x)e^x\mathds{1}_{[0,1]}(t)$, and $a(x)=1+\sin(\pi x)$.
\end{itemize}
The direct problem \eqref{eq1} is equipped with a zero Dirichlet boundary condition. The following five
weights are tested: $\mu_1=\mathds{1}_{[0.2,0.8]}$, $\mu_2=\mathds{1}_{[0.2,0.6]}$, $\mu_3=\mathds{1}_{[0.2,0.4]}$, $\mu_4=
\mathds{1}_{[0.4,0.8]}$, and $\mu_5=\mathds{1}_{[0.4,0.8]}$. {The weights share either upper or lower
bounds to facilitate the comparison of the asymptotics.}
\end{example}

\begin{figure}[hbt!]
  \begin{tabular}{ccc}
   \includegraphics[width=.45\textwidth]{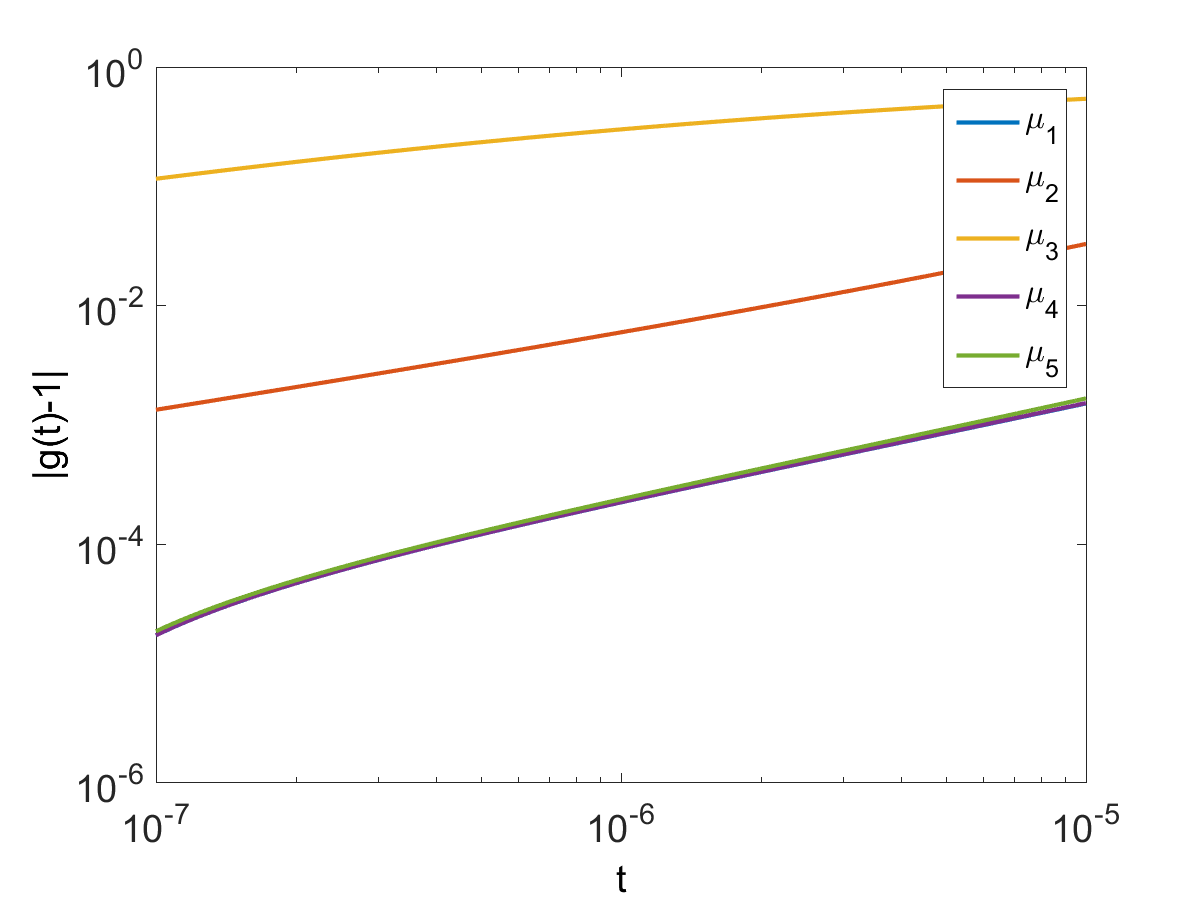} & \includegraphics[width=.45\textwidth]{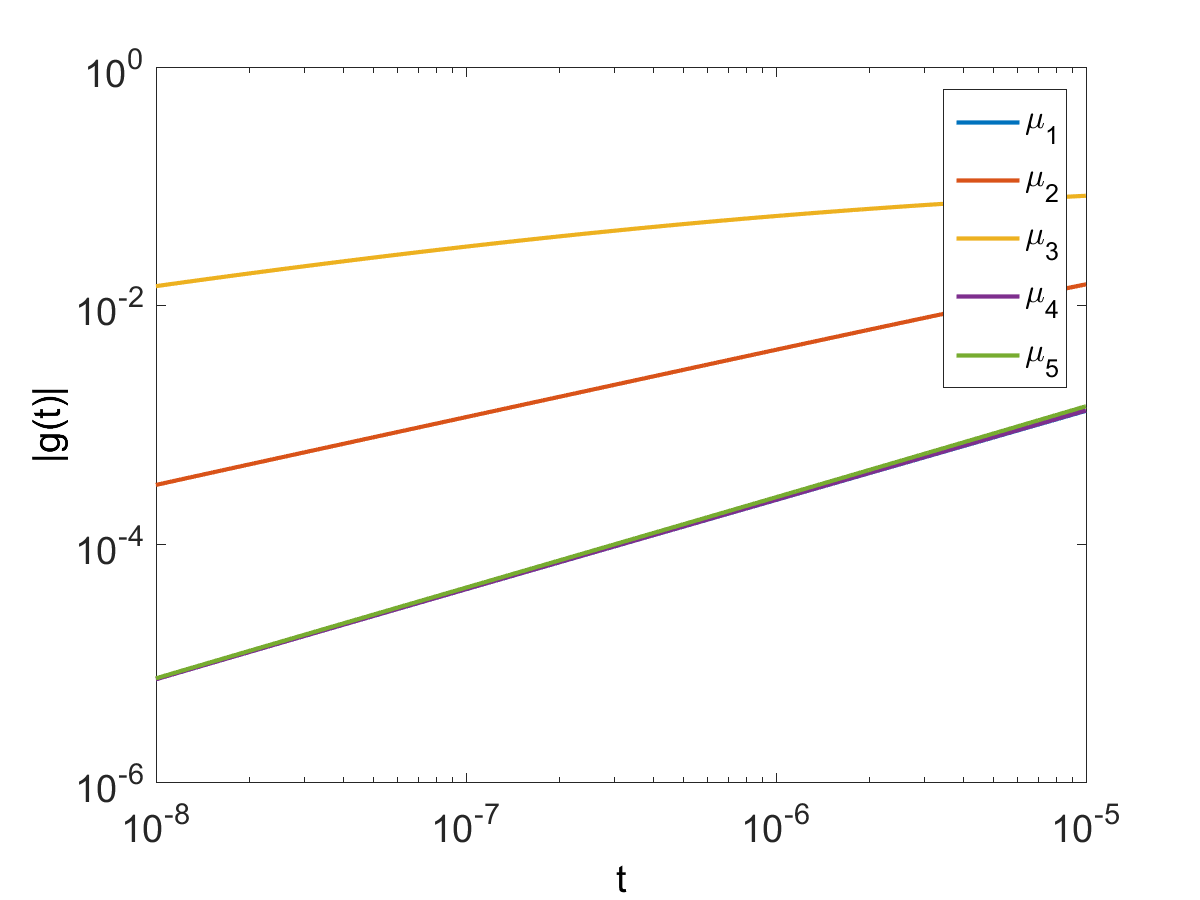}\\
   \includegraphics[width=.45\textwidth]{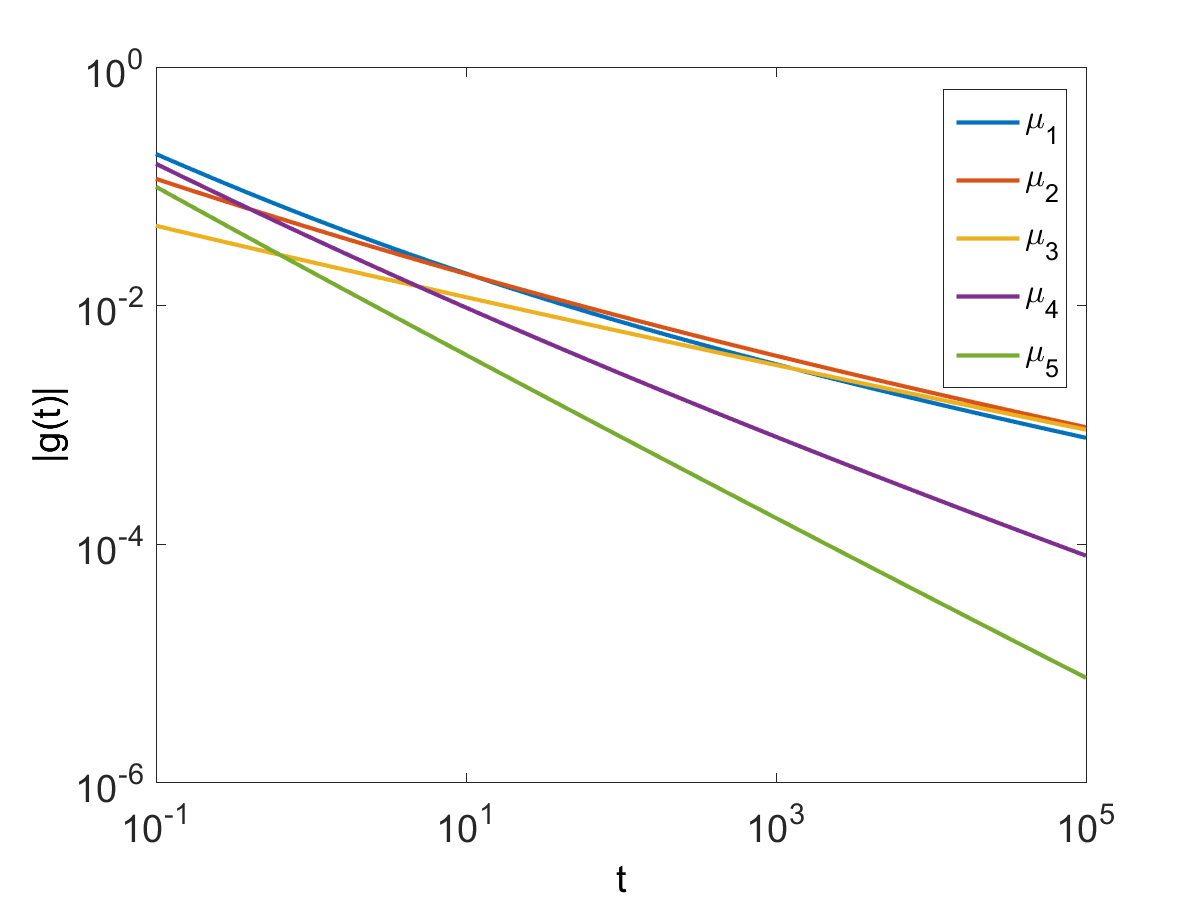} & \includegraphics[width=.45\textwidth]{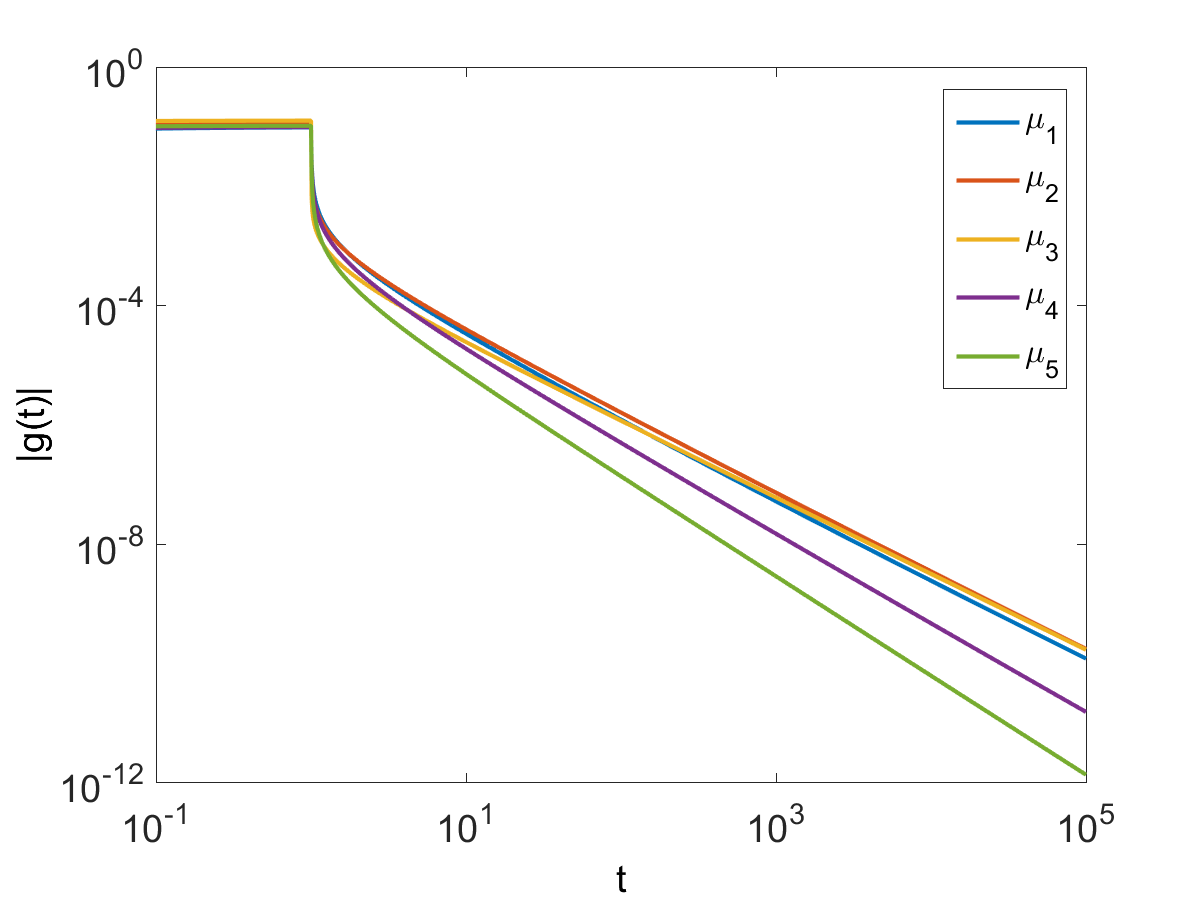}\\
   case (i) & case (ii)
  \end{tabular}
  \caption{The small (top) and large (bottom) time asymptotics of the data $g(t)=\partial_\nu u(x_0,t)$
  for Example \ref{exam:bound}, at $x_0=0$.}\label{fig:bdd}
\end{figure}

The measured data $g(t)=\partial_{\nu_a} u(x_0,t)$, at $x_0=0$, are shown in Fig. \ref{fig:bdd}. For
small time, we plot $|g(t)-g(0)|$ in order to
show the asymptotics more clearly. In case (i), the asymptotics of
$|g(t)-g(0)|$ (noting $g(0)=1$) is identical for the weights $\mu_1$, $\mu_4$ and $\mu_5$
but they differ markedly from that for $\mu_2$ and
$\mu_3$. This agrees well with the analysis: the upper bound $b_2$ of the weight support determines the small
time asymptotics. Indeed, the upper bounds of $\mu_1$, $\mu_4$ and $\mu_5$ are identical (all
equal to 0.8), whereas that of $\mu_2$ and $\mu_3$ are different. A similar behavior
is observed for case (ii). In contrast, for large time, the asymptotic of $|g(t)|$
is nearly identical for the weights $\mu_1,$ $\mu_2$ and $\mu_3$, which differs from
that for $\mu_4$ and $\mu_5$. This again agrees with the analysis: the lower bound
$b_1$ determines the large time asymptotics. Indeed, the lower bounds of $\mu_1$, $\mu_2$ and $\mu_3$
are identical (all equal to 0.2), whereas that of $\mu_4$ and $\mu_5$ are different. This observation
is the same for case (ii), except the small transition near $t=0$. In sum, the results
confirm the asymptotic behavior.

Next we recover the support bound. This is more challenging,
similar to that for multi-term case \cite{JinKian:2021pa}.
By the asymptotics in Lemma
\ref{l8}, we may recover the bound $b$ by minimizing
\begin{equation*}
  J(b,c_0,c_1) = \|c_0+c_1t^{-b} - h(x_0,t)\|^2_{L^2(t_1,t_2)},
\end{equation*}
with small $(t_1,t_2)$ for the upper bound $b_2$, and
large $t_1,t_2$ for the lower bound $b_1$. The numerical results
for Example \ref{exam:bound}(ii) are given in Table \ref{tab:bound}, where the lower
and upper bounds are estimated from the observation over the intervals $(\text{1e-6},
\text{1e-5})$ and $(\text{1e4},\text{1e5})$, respectively. Note that the estimated
lower bound is slightly more accurate than the estimated upper bound. The results for case (i)
are not shown, since the upper bound estimates are very inaccurate, similar to what
observed in the multi-term case \cite{JinKian:2021pa}. This failure might be
attributed to the fact that the interval $(\text{1e-6},\text{1e-5})$
is still too large. This can also be seen from Fig. \ref{fig:bdd}, where
the curves nearly parallel the $x$-axis, indicating a quasi-steady state. Nonetheless,
the lower bound estimates are largely the same with that for case (ii).

\begin{table}[hbt!]
  \centering
  \caption{The upper and lower bounds estimated from $g(t)$ for Example \ref{exam:bound}(ii).}\label{tab:bound}
  \begin{tabular}{c|ccccc}
   \toprule
      & $\mu_1$ & $\mu_2$ & $\mu_3$ & $\mu_4$ & $\mu_5$\\
   \midrule
  $b_1$ & 0.30 (0.20) & 0.30 (0.20) & 0.27 (0.20) &  0.48 (0.40) & 0.66 (0.60)\\
  $b_2$ & 0.75 (0.80) & 0.55 (0.60) & 0.16 (0.40) &  0.75 (0.80) & 0.76 (0.80)\\
   \bottomrule
  \end{tabular}
\end{table}

\subsection{Weight recovery}
Now we recover the weight $\mu$ from the observation $\mathcal{R}^*u(\mu)(x_0,t)$, $0<t<T$,
where $u(\mu)$ denotes the solution to problem \eqref{eq1}, when the medium is known. We
employ an iterative method to approximately minimize
\begin{equation}\label{eqn:Tikh}
  J(\mu) :=\tfrac12\|\mathcal{R}^*u(\mu)(x_0,\cdot)-g^\delta\|_{L^2(0,T)}^2.
\end{equation}
Here one may
add a penalty in order to promote distinct features, e.g., total variation for discontinuous
weight, which however is not investigated below. To minimize $J(\mu)$,
we employ the conjugate gradient method \cite{AlifanovArtyukhin:1995}, which
is a regularization method when equipped with early stopping. To apply the method,
one needs the gradient $J'(\mu)$, which can be evaluated efficiently using the adjoint technique
\cite{Cea:1986}. Assuming $\mathcal{R}$ being the Neumann trace, let $v$ be the solution
to the following adjoint problem
\begin{equation}\label{eqn:fde-adj}
  \left\{\begin{aligned}
    _t\kern-.5em^R\kern-.2em\partial_T^{[\mu]} v + \mathcal{A}v &= 0,\quad \mbox{in } \Omega\times (0,T),\\
    \partial_{\nu_a} v(x,t) & = (u(\mu)-g^\delta)\delta_{x_0}, \quad \mbox{on }\partial\Omega\times(0,T),\\
    {_tI_T^{[\mu]}}v(x,T) & = 0,\quad \mbox{in }\Omega,
  \end{aligned}  \right.
\end{equation}
where the notation $_tI_T^{[\mu]}v$ and $_t\kern-.5em^R\kern-.2em\partial_T^{[\mu]} v$
denote the right-sided Riemann-Liouville fractional integral and derivative of $v$ (with
respect to $\mu$), defined by $_tI_T^{[\mu]} v(t) = \int_t^TK_\mu(s-t)v(s)\d s$ and
$_t\kern-.5em^R\kern-.2em\partial_T^{[\mu]} v(t) = -\frac{\d}{\d t}{_tI_T^{[\mu]}}u(s)$.
Note that the Neumann data in \eqref{eqn:fde-adj}
involves a Dirac delta function $\delta_{x_0}(x)$ for some $x_0\in \partial\Omega$ and it
should be interpreted in the sense of distribution. When implementing the Galerkin FEM, one should
evaluate the integral on the boundary carefully, which results in a
term supported at $x_0$ only, cf. \eqref{eqn:fde-adj-weak} below. By changing variables and using the zero terminal data,
problem \eqref{eqn:fde-adj} can be transformed into an initial value problem with a
Djrbashian-Caputo type derivative and zero initial data, which can then be solved using the standard
scheme in the appendix.

We have the following representation of the gradient $J'(\mu)$.
\begin{proposition}
The gradient $J'(\mu)$ of the functional $J(\mu)$ is given by
\begin{equation}\label{eqn:grad}
  J'(\mu):(0,1)\ni\alpha\mapsto -\int_0^T\int_\Omega v \partial_t^\alpha u \d x \d t,
\end{equation}
where $v$ is the solution to the adjoint problem \eqref{eqn:fde-adj} and, for all $\alpha\in(0,1)$,  $\partial_t^\alpha u$ is defined by
$\partial_t^\alpha u(x,t):=\partial_t\int_0^t\frac{(t-s)^{-\alpha}}{\Gamma(1-\alpha)}(u(x,s)-u_0(x))\d s$, for $(x,t)\in\Omega\times(0,T).$
\end{proposition}
\begin{proof}
The derivation follows by a standard procedure. Let $u\equiv u(\mu)$. By definition, the directional derivative $J'(\mu)[h]$
of the functional $J$ with respect to the weight $\mu$ in the direction $h\in L^2(0,1)$ is given by
\begin{equation*}
  J'(\mu)[h] = (u_h(x_0,t),u(x_0,t)-g^\delta(t))_{L^2(0,T)},
\end{equation*}
where $u_h$ is the derivative of $u$ with respect to the weight $\mu$ in the direction $h$. Clearly $u_h$ satisfies
\begin{equation}\label{eqn:fde-lin0}
  \left\{\begin{aligned}
    \partial_t^{[\mu]} u_h + \mathcal{A} u_h &= -\int_0^1 h(\alpha)\partial_t^\alpha u\d\alpha, \quad \mbox{in }\Omega\times (0,T),\\
    \partial_{\nu_a} u_h(x,t) & = 0, \quad \mbox{on }\partial\Omega\times(0,T),\\
    u_h(x,0) & = 0,\quad \mbox{in }\Omega.
  \end{aligned}  \right.
\end{equation}
Multiplying \eqref{eqn:fde-lin0} with a test function $\phi(x,t)$, integrating over $\Omega\times (0,T)$
and then applying integration by parts yield
\begin{equation}\label{eqn:fde-lin}
  \int_0^T\!\!\int_\Omega\!\! \Big(\int_0^1\partial_t^\alpha u_h\mu(\alpha)\d \alpha\Big)\phi+ a\nabla u_h\cdot\nabla \phi \d x\d t = \int_0^T\!\!\int_\Omega \Big(-\int_0^1\partial_t^\alpha u h(\alpha)\d \alpha\Big)\phi\d x \d t.
\end{equation}
Meanwhile, the weak formulation for the adjoint solution $v$ is given by
\begin{equation}\label{eqn:fde-adj-weak}
  \int_0^T\!\!\int_\Omega\!\! \Big(\int_0^1 {_t\kern-.5em^R\kern-.2em\partial_T^\alpha}v\mu(\alpha)\d\alpha\Big)\phi + a\nabla v\cdot\nabla \phi \d x\d t = \int_0^T(u(\mu)(x_0,t)-g^\delta)\phi(x_0,t)\d t,
\end{equation}
where we recall that for all $\alpha\in(0,1)$,  ${_t\kern-.5em^R\kern-.2em\partial_T^\alpha}v$ is defined by
$${_t\kern-.5em^R\kern-.2em\partial_T^\alpha}v(x,t):=-\partial_t\left(\int_t^T\frac{(t-s)^{-\alpha}}{\Gamma(1-\alpha)}v(x,s)\d s\right),\quad (x,t)\in\Omega\times(0,T).$$
Since $u_h(0)=0$ and $v(T)=0$, the following identity holds (\cite[p. 76, Lemma 2.7]{KilbasSrivastavaTrujillo:2006} or \cite[Lemma 2.6]{Jin:2021book})
$$\int_0^T\int_\Omega v\partial_t^\alpha u_h \d x\d t = \int_0^T\int_\Omega u_h{_t\kern-.5em^R\kern-.2em\partial_T^\alpha} v\d x\d t,\quad \alpha\in(0,1).$$
Using this identity, taking $\phi=v$ in \eqref{eqn:fde-lin} and $\phi=u_h$ in \eqref{eqn:fde-adj-weak}
and subtracting the resulting identities give
\begin{equation*}
  \int_0^1 \Big(-\int_0^T\int_\Omega \partial_t^\alpha u v\d x \d t \Big)h(\alpha)\d\alpha = \int_0^Tu_h(x_0,t) (u(\mu)(x_0,t)-g^\delta)\d t.
\end{equation*}
This and the definition of the derivative $J'(\mu)$ show the desired assertion.
\end{proof}

Now we can describe the conjugate gradient method \cite{AlifanovArtyukhin:1995} for minimizing the functional $J$;
see Algorithm \ref{alg:cgm} for the complete procedure. In the algorithm, Steps 6-7 compute the conjugate descent
direction, and Step 9 computes the optimal step size using the sensitivity problem. The operator $P_+$ at step 10
gives the projection into the positive quadrant. For each update, the algorithm involves three forward solves
(direct, adjoint and sensitivity problems), which represent the main computational effort. For the stopping criterion
at Step 11, one may employ the discrepancy principle \cite{Morozov:1966,ItoJin:2015}. Note that this algorithm has
been adopted in \cite{LiuSunYamamoto:2021} recently, where the gradient is computed using a finite difference
approximation, which is generally less efficient than the adjoint technique. In the experiments, we smooth the
gradient $J'(\mu)$ by the negative Dirichlet Laplacian, which is analogous to the $H^1(0,1)$ penalty in
\cite{LiuSunYamamoto:2021}. Accordingly, at Step 6 of the algorithm, we compute the conjugate
coefficient $\gamma$ in the $H^1(0,1)$-seminorm, where $\partial_\alpha$ denotes taking the first-order derivative in $\alpha$.
This approach is suitable for recovering a smooth weight.

\begin{algorithm}[hbt!]
  \caption{Conjugate gradient method for minimizing the functional $J$.\label{alg:cgm}}
  \begin{algorithmic}[1]
    \STATE Initialize $\mu^0$, and set $k=0$.
    \FOR{$k=0,\ldots,K$}
    \STATE Solve for $u^k$ from problem \eqref{eq1} with $\mu=\mu^k$.
    \STATE Solve for $v^k$ from problem \eqref{eqn:fde-adj} with $r^k=u^k(x_0,t)-g^\delta$.
    \STATE Compute the gradient $J'(\mu^k)$  by \eqref{eqn:grad}.
    \STATE Compute the conjugate coefficient $\gamma^k$ by $\gamma^0=0$ and
    \begin{equation*}
      \gamma^k = \frac{\|\partial_\alpha J'(\mu^k)\|_{L^2(0,1)}^2}{\|\partial_\alpha J'(\mu^{k-1})\|_{L^2(0,1)}^2}
      , \quad k\geq 1.
    \end{equation*}
    \STATE Compute the conjugate direction $d^k$ by $d^k=-J'(\mu^k)+\gamma^kd^{k-1}$.
    \STATE Solve for $u_{d^k}$ from problem \eqref{eqn:fde-lin0} with $h=d^k$.
    \STATE Compute the step size $s^k$ by
    \begin{equation*}
      s^k = -\frac{(u_{d^k}(x_0,\cdot),r^k)_{L^2(0,T)}}{\|u_{d^k}(x_0,\cdot)\|^2_{L^2(0,T)}}.
    \end{equation*}
    \STATE Update the source component $\mu^{k+1}=P_+(\mu^k+s^k d^k)$.
    \STATE Check the stopping criterion.
    \ENDFOR
  \end{algorithmic}
\end{algorithm}

The next example illustrates the weight recovery.
\begin{example}\label{exam:weight}
In this example, we take the domain $\Omega=(0,1)$, terminal time $T= 1$,  diffusion coefficient
$a(x)=1+x(1-x)$, initial condition $u_0(x)=x(1-x)e^x$. The direct problem \eqref{eq1} is equipped with
the Neumann boundary condition, which at $x=0$ and $x=1$ are taken to be $0$ and $1$, respectively.
Consider the following two weights: $\rm(i)$ $\mu(\alpha)= \alpha(1-\alpha)^2e^{2\alpha}$ and
$\rm(ii)$ $\mu(\alpha)=2\min(\alpha,1-\alpha)$. The Dirichlet data at the point $x_0=0$ is used for
the recovery. Algorithm \ref{alg:cgm} is initialized with $\mu^0=\frac{1}{100}\sin(\pi\alpha)$.
\end{example}

\begin{table}[hbt!]
  \centering
  \caption{The $L^2(0,1)$ error of the recovered weight $\hat \mu$ for Example \ref{exam:weight} at different noise levels.\label{tab:weight}}
  \setlength{\tabcolsep}{3pt}
  \begin{tabular}{ccccccc}
  \toprule
    $\epsilon$  & 0 & 1e-3 & 3e-3 & 1e-2 & 3e-2 & 5e-2\\
  \midrule
  (i)  &  2.50e-3(22) & 4.09e-3(22) &  8.82e-3(50) &   1.62e-2(20) &  4.16e-3(17) &  9.33e-3(17) \\
  (ii) &  4.49e-2(26) & 4.48e-2(26) &  4.63e-2(20) &   5.58e-2(18) &  4.96e-2(17) &  5.99e-2(14)\\
  \bottomrule
  \end{tabular}
\end{table}

\begin{figure}
  \centering
  \setlength{\tabcolsep}{0pt}
  \begin{tabular}{ccc}
   \includegraphics[width=.32\textwidth]{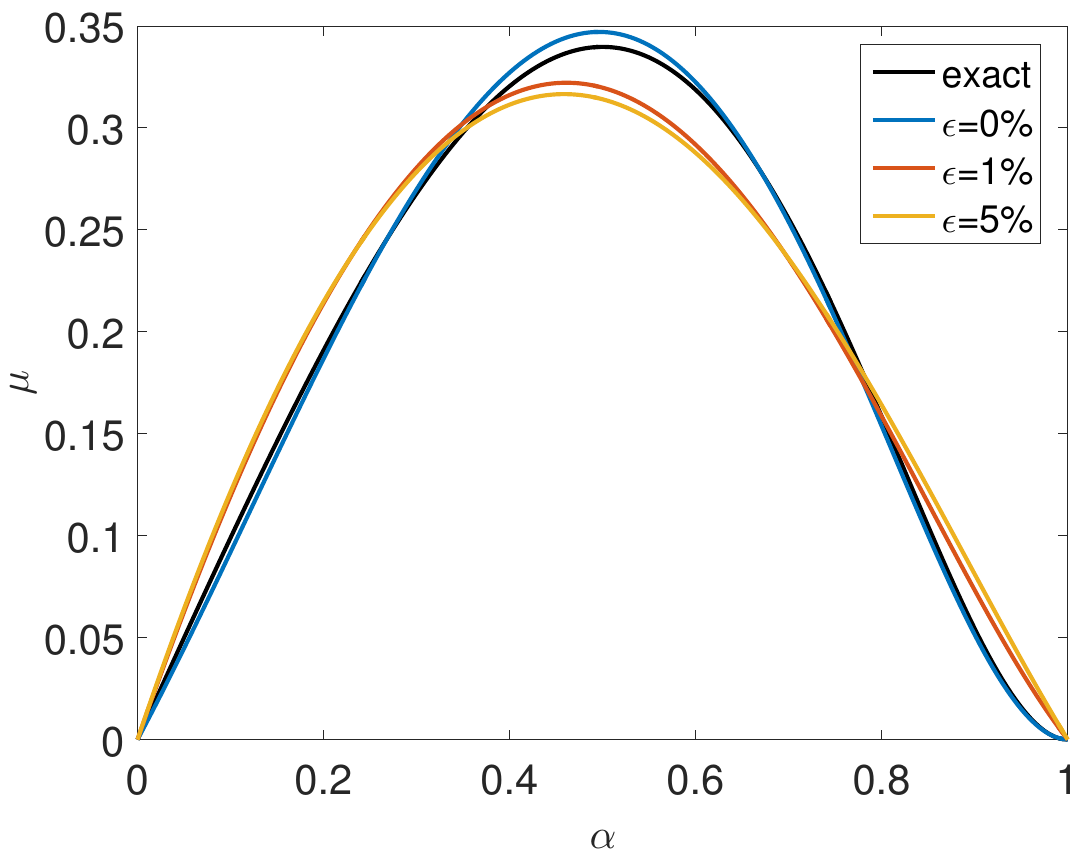} & \includegraphics[width=.32\textwidth]{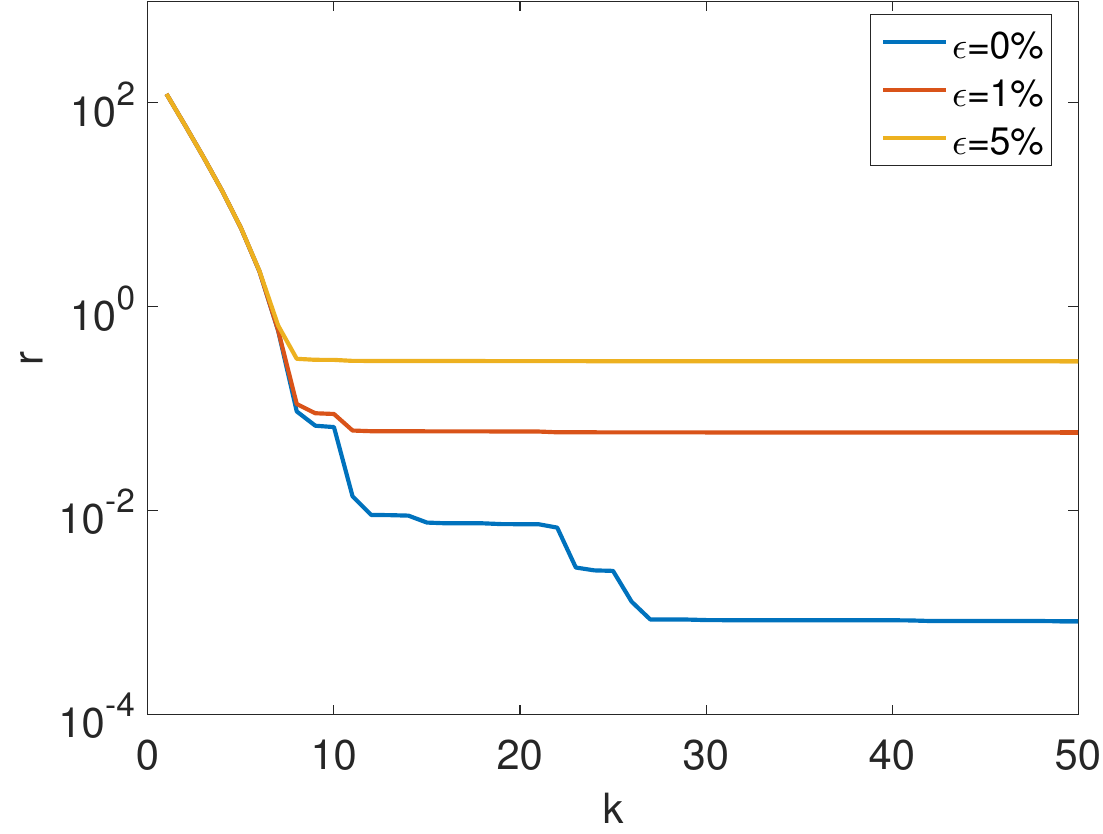} & \includegraphics[width=.32\textwidth]{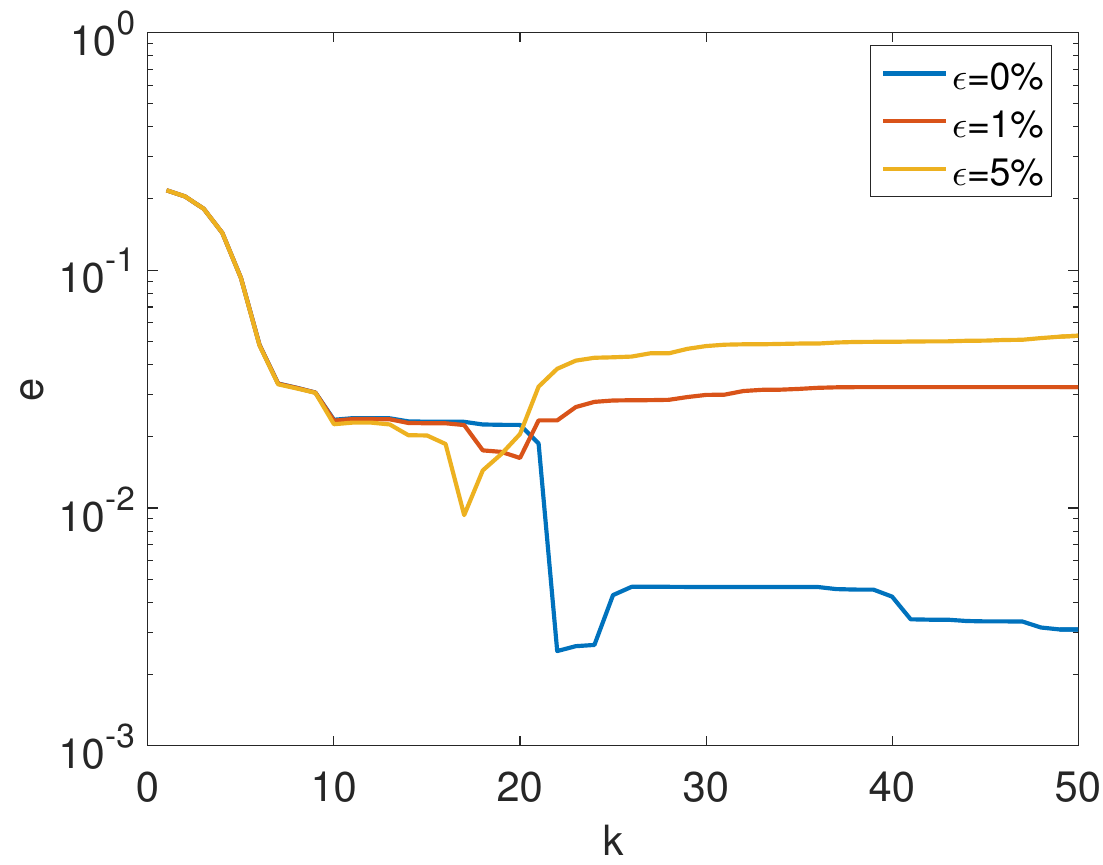}\\
   (a) reconstruction & (b) residual & (c) $L^2(0,1)$ error
  \end{tabular}
  \caption{Numerical results for Example \ref{exam:weight}(i) by the conjugate gradient method.}\label{fig:weight-smooth}
\end{figure}

\begin{figure}
  \centering
  \setlength{\tabcolsep}{0pt}
  \begin{tabular}{ccc}
   \includegraphics[width=.32\textwidth]{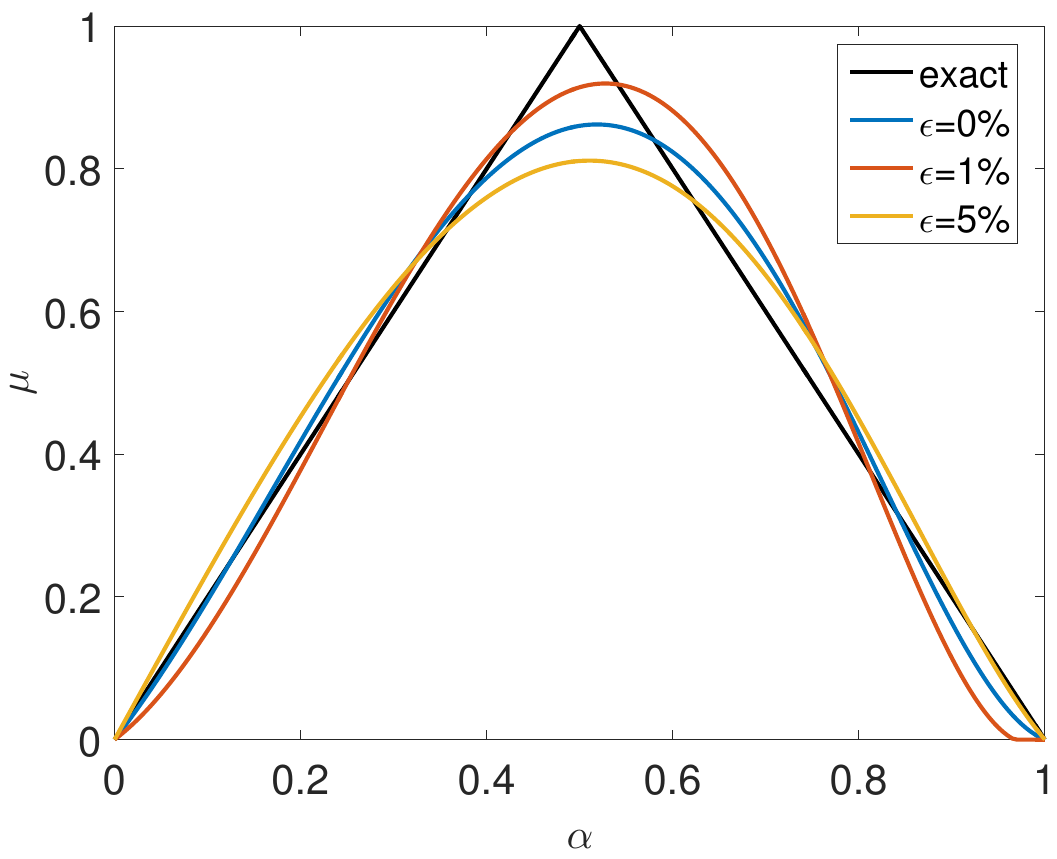} & \includegraphics[width=.32\textwidth]{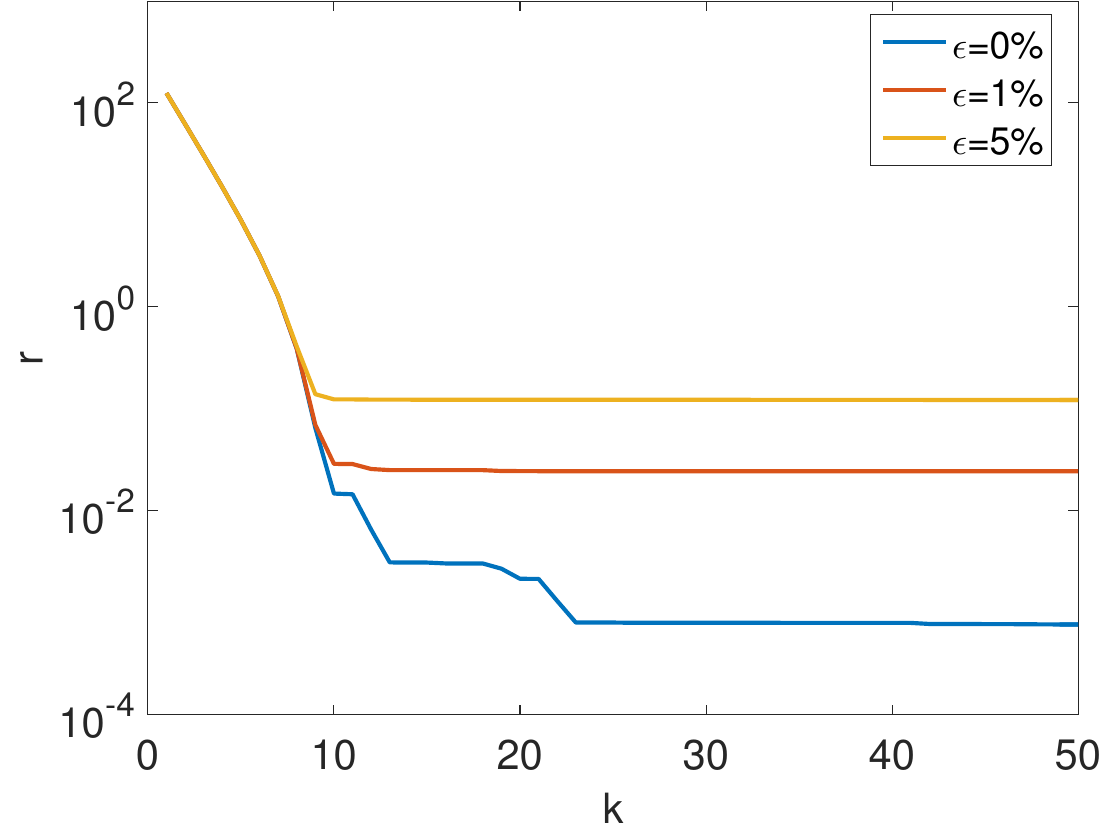} & \includegraphics[width=.32\textwidth]{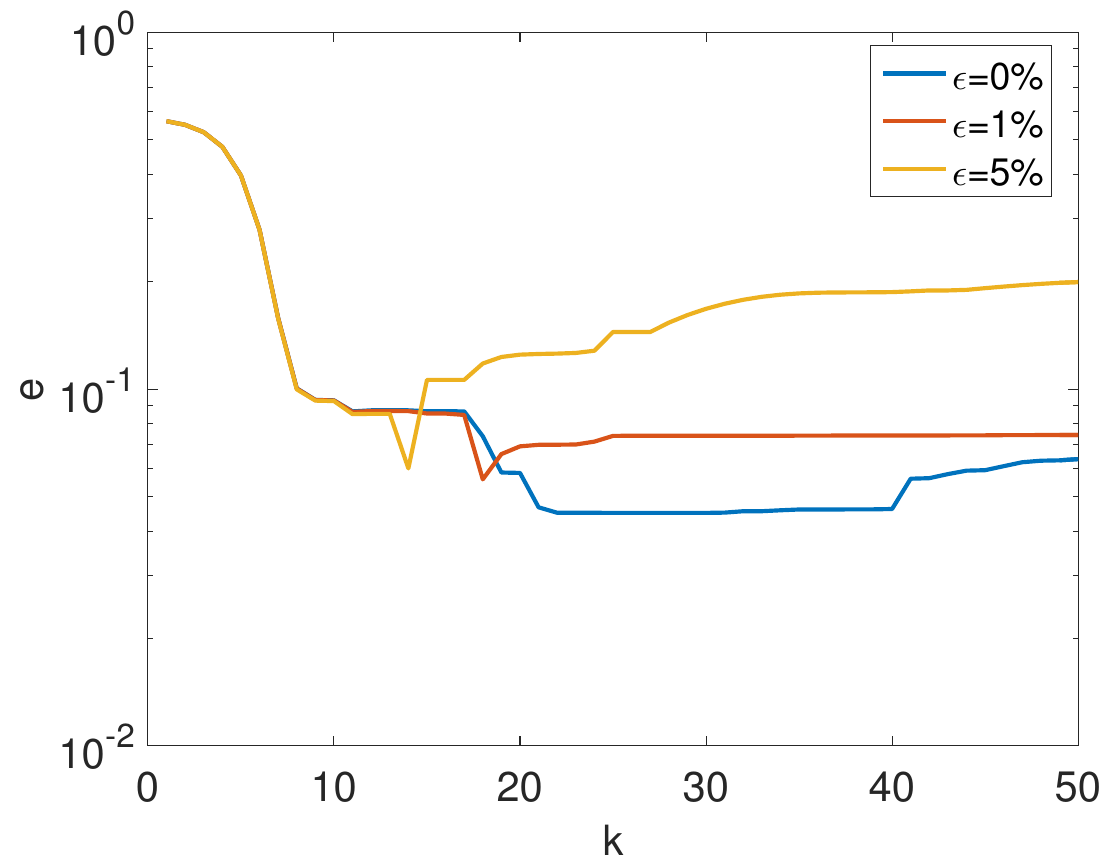}\\
   (a) reconstruction & (b) residual & (c) $L^2(0,1)$ error
  \end{tabular}
  \caption{Numerical results for Example \ref{exam:weight}(ii) by the conjugate gradient method.}\label{fig:weight-nonsmooth}
\end{figure}

The numerical results are shown in Table \ref{tab:weight} and Figs.
\ref{fig:weight-smooth} and \ref{fig:weight-nonsmooth}, where the stopping index $k^*$ (in the bracket)
is taken such that the $L^2(0,1)$ error $e=\|\hat\mu-\mu^\dag\|_{L^2(0,1)}$ of the recovered
weight $\hat \mu$ (with respect to the exact one $\mu^\dag$) is smallest. The
weight $\mu$ can be fairly reconstructed for the noise level $\epsilon$ up to $5\%$, agreeing with the results
in Table \ref{tab:weight}. The conjugate gradient method converges steadily, irrespective of
the noise level. However, early stopping is needed in order to obtain high-quality reconstructions,
due to the typical semiconvergence behaviour of the method: the error first decreases steadily
and then it increases as the iteration further proceeds. This observation holds for both smooth and nonsmooth
weights, which shows the feasibility of recovering the weight $\mu$ from the observation at one
single point, when all other problem data are fully known, confirming the assertion of Theorem \ref{thm:known}.

\section{Conclusion}\label{sec:concl}

In this work we have investigated the inverse problem of recovering information about the
weight in a distributed-order time-fractional diffusion model from one-point observation on the
boundary. We prove the support bound recovery when the medium is unknown, and the weight
recovery when the medium is known. The latter result complements an existing uniqueness result
for a known medium. Further
we have also presented numerical experiments to complement the analysis, partly confirming the
unique recovery of either support boundary or weight. The experiments show
the significant computational challenge with accurately recovering support bounds and weights,
especially the upper bound of the support. Theoretically, it is of much
interest to study the case of a nonseparable source. This would require alternative analysis
techniques than Titchmarsh convolution theorem, which has played a crucial role in the analysis of this work.

\appendix
\section{Numerical scheme for the direct problem \eqref{eq1}}
Now we describe a numerical scheme for approximating the direct problem \eqref{eq1} for
completeness; see the book \cite{JinZhou:2022book} for further details on the numerical approximation of
the standard time-fractional diffusion model. For the spatial discretization,
we employ the standard Galerkin finite element method with the continuous piecewise linear finite element space
$X_h$, subordinated to a shape regular
triangulation of the domain $\Omega$. Let $M_h$ and $S_h$ be the corresponding mass and stiffness
matrices, respectively. For any $N\in\mathbb{N}$ total number of time steps, let $\tau=N^{-1}T$ be the
time step size, and $t_n=n\tau$, $n=0,\ldots, N$, the time grid. The classical L1 scheme
(cf. \cite{LinXu:2007} or \cite[Section 4.1]{JinZhou:2022book}) approximates
the Djrbashian-Caputo fractional derivative $\partial_t^\alpha u(t_n)$, $n\geq1$, by
\begin{align*}
  \partial_t^\alpha u^n     & \approx \sum_{j=0}^n b_{j,n}^{(\alpha)}u^{n-j},
\end{align*}
with the weights $b_{j,n}^{(\alpha)}$ given by
\begin{equation*}
  b_{j,n}^{(\alpha)} = \frac{1}{\Gamma(2-\alpha)\tau^\alpha}\left\{\begin{aligned}
    1, & \quad j =0,\\
    (j+1)^{1-\alpha}+(j-1)^{1-\alpha}-2j^{1-\alpha}, &\quad j =1,\ldots,n-1,\\
    n^{1-\alpha}-(n-1)^{1-\alpha}, &\quad j = n.
  \end{aligned}\right.
\end{equation*}
For the distributed-order fractional derivative $\partial_t^{[\mu]} u^n$, one may use quadrature
to approximate the integral with respect to the weight $\mu(\alpha)$. For example, the trapezoidal rule gives
\begin{align*}
  \partial_t^{[\mu]} u^n &= \int_0^1\partial_t^\alpha u(t_n) \mu(\alpha)\d \alpha \approx
   \delta_\alpha\sum_{i=0}^{N_\alpha}c_i\mu(\alpha_i) \partial_t^{\alpha_i} u^n,
\end{align*}
with $c_0=c_{N_\alpha}=1/2$, and $c_i=1$ for $i=1,\ldots,N_\alpha-1$, and $\delta_\alpha =N_\alpha^{-1}$.
These two approximations lead to
\begin{align*}
 \partial_t^{[\mu]} u^n  \approx \sum_{j=0}^{n}p_{n-j}^{[\mu]}(u^j-u^0), \quad \mbox{with } p_j^{[\mu]} =\delta_\alpha \sum_{i=0}^{N_\alpha} c_i\mu(\alpha_i) b_{j,n}^{(\alpha_i)}.
\end{align*}
Note that the complexity of the scheme is comparable with the standard one, except computing the
weights $p_j^{[\mu]}$. The fully discrete scheme reads: with $U_h^0=P_hu_0$ and $F_h^n=P_hf(t_n)$ (with
$P_h$ being the $L^2(\Omega)$ projection on $X_h$), find $U_h^i\in X_h$ such that for $n=1,2,\ldots,N$,
\begin{align*}
   (p_0^{[\mu]}M_h+S_h)U_h^n = p_0^{[\mu]}U^0_h-\sum_{j=1}^{n-1}p_j^{[\mu]}(U_h^{n-j}-U_h^0) + F_h^n.
\end{align*}
Note that the computational complexity per time step of the scheme increases linearly with the
time step $n$, which results in a quadratic complexity of the overall scheme. For large time
simulation, one may use the sum of exponential approximation to reduce the complexity
\cite[Section 4.2]{JinZhou:2022book}.

\bibliographystyle{abbrv}
\bibliography{frac}
\end{document}